\documentclass[11pt]{article}
\usepackage{amsmath, upgreek}
\usepackage{amssymb}
\usepackage{color}
\usepackage{amscd}
\usepackage{xspace}
\usepackage{verbatim}
\usepackage{graphicx}
\renewcommand{\theequation}{\thesection.\arabic{equation}
}

\title{\bf Weak solutions of  semilinear elliptic  equations
with\\  Leray-Hardy potential
and measure data  }
\author{{\bf Huyuan Chen\footnote{\noindent Department of Mathematics, Jiangxi Normal University,
Nanchang 330022, China. E-mail: chenhuyuan@yeah.net}} \\[4mm]
 {\bf Laurent V\'eron \footnote{\noindent
Laboratoire de Math\'{e}matiques et Physique Th\'{e}orique, Universit\'e de Tours, 37200 Tours, France. E-mail: veronl@univ-tours.fr}}
}

\date{}
\begin{document}
 \maketitle


\newcommand{\txt}[1]{\;\text{ #1 }\;}
\newcommand{\tbf}{\textbf}
\newcommand{\tit}{\textit}
\newcommand{\tsc}{\textsc}
\newcommand{\trm}{\textrm}
\newcommand{\mbf}{\mathbf}
\newcommand{\mrm}{\mathrm}
\newcommand{\bsym}{\boldsymbol}
\newcommand{\scs}{\scriptstyle}
\newcommand{\sss}{\scriptscriptstyle}
\newcommand{\txts}{\textstyle}
\newcommand{\dsps}{\displaystyle}
\newcommand{\fnz}{\footnotesize}
\newcommand{\scz}{\scriptsize}
\newcommand{\be}{\begin{equation}}
\newcommand{\bel}[1]{\begin{equation}\label{#1}}
\newcommand{\ee}{\end{equation}}
\newcommand{\eqnl}[2]{\begin{equation}\label{#1}{#2}\end{equation}}
\newcommand{\barr}{\begin{eqnarray}}
\newcommand{\earr}{\end{eqnarray}}
\newcommand{\bars}{\begin{eqnarray*}}
\newcommand{\ears}{\end{eqnarray*}}
\newcommand{\nnu}{\nonumber \\}
\newtheorem{subn}{\name}
\renewcommand{\thesubn}{}
\newcommand{\bsn}[1]{\def\name{#1}\begin{subn}}
\newcommand{\esn}{\end{subn}}
\newtheorem{sub}{\name}[section]
\newcommand{\dn}[1]{\def\name{#1}}   
\newcommand{\bs}{\begin{sub}}
\newcommand{\es}{\end{sub}}
\newcommand{\bsl}[1]{\begin{sub}\label{#1}}
\newcommand{\bth}[1]{\def\name{Theorem}
\begin{sub}\label{t:#1}}
\newcommand{\blemma}[1]{\def\name{Lemma}
\begin{sub}\label{l:#1}}
\newcommand{\bcor}[1]{\def\name{Corollary}
\begin{sub}\label{c:#1}}
\newcommand{\bdef}[1]{\def\name{Definition}
\begin{sub}\label{d:#1}}
\newcommand{\bprop}[1]{\def\name{Proposition}
\begin{sub}\label{p:#1}}

\newcommand{\aand}{\quad\mbox{and}\quad}
\newcommand{\M}{{\cal M}}
\newcommand{\A}{{\cal A}}
\newcommand{\B}{{\cal B}}
\newcommand{\I}{{\cal I}}
\newcommand{\J}{{\cal J}}
\newcommand{\D}{\displaystyle}
\newcommand{\RR}{ I\!\!R}
\newcommand{\C}{\mathbb{C}}
\newcommand{\R}{\mathbb{R}}
\newcommand{\Z}{\mathbb{Z}}
\newcommand{\N}{\mathbb{N}}
\newcommand{\T}{{\rm T}^n}
\newcommand{\cuad}{{\sqcap\kern-.68em\sqcup}}
\newcommand{\abs}[1]{\mid #1 \mid}
\newcommand{\norm}[1]{\|#1\|}
\newcommand{\equ}[1]{(\ref{#1})}
\newcommand\rn{\mathbb{R}^N}
\renewcommand{\theequation}{\thesection.\arabic{equation}}
\newtheorem{definition}{Definition}[section]
\newtheorem{theorem}{Theorem}[section]
\newtheorem{proposition}{Proposition}[section]
\newtheorem{example}{Example}[section]
\newtheorem{proof}{proof}[section]
\newtheorem{lemma}{Lemma}[section]
\newtheorem{corollary}{Corollary}[section]
\newtheorem{remark}{Remark}[section]
\newcommand{\bremark}{\begin{remark} \em}
\newcommand{\eremark}{\end{remark} }
\newtheorem{claim}{Claim}


\newcommand{\rth}[1]{Theorem~\ref{t:#1}}
\newcommand{\rlemma}[1]{Lemma~\ref{l:#1}}
\newcommand{\rcor}[1]{Corollary~\ref{c:#1}}
\newcommand{\rdef}[1]{Definition~\ref{d:#1}}
\newcommand{\rprop}[1]{Proposition~\ref{p:#1}}
\newcommand{\BA}{\begin{array}}
\newcommand{\EA}{\end{array}}
\newcommand{\BAN}{\renewcommand{\arraystretch}{1.2}
\setlength{\arraycolsep}{2pt}\begin{array}}
\newcommand{\BAV}[2]{\renewcommand{\arraystretch}{#1}
\setlength{\arraycolsep}{#2}\begin{array}}
\newcommand{\BSA}{\begin{subarray}}
\newcommand{\ESA}{\end{subarray}}
\newcommand{\BAL}{\begin{aligned}}
\newcommand{\EAL}{\end{aligned}}
\newcommand{\BALG}{\begin{alignat}}
\newcommand{\EALG}{\end{alignat}}
\newcommand{\BALGN}{\begin{alignat*}}
\newcommand{\EALGN}{\end{alignat*}}
\newcommand{\note}[1]{\textit{#1.}\hspace{2mm}}
\newcommand{\Proof}{\note{Proof}}
\newcommand{\qeda}{\hspace{10mm}\hfill $\square$}
\newcommand{\qed}{\\
${}$ \hfill $\square$}
\newcommand{\Remark}{\note{Remark}}
\newcommand{\modin}{$\,$\\[-4mm] \indent}
\newcommand{\forevery}{\quad \forall}
\newcommand{\set}[1]{\{#1\}}
\newcommand{\setdef}[2]{\{\,#1:\,#2\,\}}
\newcommand{\setm}[2]{\{\,#1\mid #2\,\}}
\newcommand{\mt}{\mapsto}
\newcommand{\lra}{\longrightarrow}
\newcommand{\lla}{\longleftarrow}
\newcommand{\llra}{\longleftrightarrow}
\newcommand{\Lra}{\Longrightarrow}
\newcommand{\warrow}{\rightharpoonup}
\newcommand{
\paran}[1]{\left (#1 \right )}
\newcommand{\sqbr}[1]{\left [#1 \right ]}
\newcommand{\curlybr}[1]{\left \{#1 \right \}}
\newcommand{
\paranb}[1]{\big (#1 \big )}
\newcommand{\lsqbrb}[1]{\big [#1 \big ]}
\newcommand{\lcurlybrb}[1]{\big \{#1 \big \}}
\newcommand{\absb}[1]{\big |#1\big |}
\newcommand{\normb}[1]{\big \|#1\big \|}
\newcommand{
\paranB}[1]{\Big (#1 \Big )}
\newcommand{\absB}[1]{\Big |#1\Big |}
\newcommand{\normB}[1]{\Big \|#1\Big \|}
\newcommand{\produal}[1]{\langle #1 \rangle}

\newcommand{\thkl}{\rule[-.5mm]{.3mm}{3mm}}
\newcommand{\thknorm}[1]{\thkl #1 \thkl\,}
\newcommand{\trinorm}[1]{|\!|\!| #1 |\!|\!|\,}
\newcommand{\bang}[1]{\langle #1 \rangle}
\def\angb<#1>{\langle #1 \rangle}
\newcommand{\vstrut}[1]{\rule{0mm}{#1}}
\newcommand{\rec}[1]{\frac{1}{#1}}
\newcommand{\opname}[1]{\mbox{\rm #1}\,}
\newcommand{\supp}{\opname{supp}}
\newcommand{\dist}{\opname{dist}}
\newcommand{\myfrac}[2]{{\displaystyle \frac{#1}{#2} }}
\newcommand{\myint}[2]{{\displaystyle \int_{#1}^{#2}}}
\newcommand{\mysum}[2]{{\displaystyle \sum_{#1}^{#2}}}
\newcommand {\dint}{{\displaystyle \myint\!\!\myint}}
\newcommand{\q}{\quad}
\newcommand{\qq}{\qquad}
\newcommand{\hsp}[1]{\hspace{#1mm}}
\newcommand{\vsp}[1]{\vspace{#1mm}}
\newcommand{\ity}{\infty}
\newcommand{\prt}{\partial}
\newcommand{\sms}{\setminus}
\newcommand{\ems}{\emptyset}
\newcommand{\ti}{\times}
\newcommand{\pr}{^\prime}
\newcommand{\ppr}{^{\prime\prime}}
\newcommand{\tl}{\tilde}
\newcommand{\sbs}{\subset}
\newcommand{\sbeq}{\subseteq}
\newcommand{\nind}{\noindent}
\newcommand{\ind}{\indent}
\newcommand{\ovl}{\overline}
\newcommand{\unl}{\underline}
\newcommand{\nin}{\not\in}
\newcommand{\pfrac}[2]{\genfrac{(}{)}{}{}{#1}{#2}}

\def\ga{\alpha}     \def\gb{\beta}       \def\gg{\gamma}
\def\gc{\chi}       \def\gd{\delta}      \def\ge{\epsilon}
\def\gth{\theta}                         \def\vge{\varepsilon}
\def\gf{\phi}       \def\vgf{\varphi}    \def\gh{\eta}
\def\gi{\iota}      \def\gk{\kappa}      \def\gl{\lambda}
\def\gm{\mu}        \def\gn{\nu}         \def\gp{\pi}
\def\vgp{\varpi}    \def\gr{\rho}        \def\vgr{\varrho}
\def\gs{\sigma}     \def\vgs{\varsigma}  \def\gt{\tau}
\def\gu{\upsilon}   \def\gv{\vartheta}   \def\gw{\omega}
\def\gx{\xi}        \def\gy{\psi}        \def\gz{\zeta}
\def\Gg{\Gamma}     \def\Gd{\Delta}      \def\Gf{\Phi}
\def\Gth{\Theta}
\def\Gl{\Lambda}    \def\Gs{\Sigma}      \def\Gp{\Pi}
\def\Gw{\Omega}     \def\Gx{\Xi}         \def\Gy{\Psi}

\def\CS{{\mathcal S}}   \def\CM{{\mathcal M}}   \def\CN{{\mathcal N}}
\def\CR{{\mathcal R}}   \def\CO{{\mathcal O}}   \def\CP{{\mathcal P}}
\def\CA{{\mathcal A}}   \def\CB{{\mathcal B}}   \def\CC{{\mathcal C}}
\def\CD{{\mathcal D}}   \def\CE{{\mathcal E}}   \def\CF{{\mathcal F}}
\def\CG{{\mathcal G}}   \def\CH{{\mathcal H}}   \def\CI{{\mathcal I}}
\def\CJ{{\mathcal J}}   \def\CK{{\mathcal K}}   \def\CL{{\mathcal L}}
\def\CT{{\mathcal T}}   \def\CU{{\mathcal U}}   \def\CV{{\mathcal V}}
\def\CZ{{\mathcal Z}}   \def\CX{{\mathcal X}}   \def\CY{{\mathcal Y}}
\def\CW{{\mathcal W}} \def\CQ{{\mathcal Q}}
\def\BBA {\mathbb A}   \def\BBb {\mathbb B}    \def\BBC {\mathbb C}
\def\BBD {\mathbb D}   \def\BBE {\mathbb E}    \def\BBF {\mathbb F}
\def\BBG {\mathbb G}   \def\BBH {\mathbb H}    \def\BBI {\mathbb I}
\def\BBJ {\mathbb J}   \def\BBK {\mathbb K}    \def\BBL {\mathbb L}
\def\BBM {\mathbb M}   \def\BBN {\mathbb N}    \def\BBO {\mathbb O}
\def\BBP {\mathbb P}   \def\BBR {\mathbb R}    \def\BBS {\mathbb S}
\def\BBT {\mathbb T}   \def\BBU {\mathbb U}    \def\BBV {\mathbb V}
\def\BBW {\mathbb W}   \def\BBX {\mathbb X}    \def\BBY {\mathbb Y}
\def\BBZ {\mathbb Z}

\def\GTA {\mathfrak A}   \def\GTB {\mathfrak B}    \def\GTC {\mathfrak C}
\def\GTD {\mathfrak D}   \def\GTE {\mathfrak E}    \def\GTF {\mathfrak F}
\def\GTG {\mathfrak G}   \def\GTH {\mathfrak H}    \def\GTI {\mathfrak I}
\def\GTJ {\mathfrak J}   \def\GTK {\mathfrak K}    \def\GTL {\mathfrak L}
\def\GTM {\mathfrak M}   \def\GTN {\mathfrak N}    \def\GTO {\mathfrak O}
\def\GTP {\mathfrak P}   \def\GTR {\mathfrak R}    \def\GTS {\mathfrak S}
\def\GTT {\mathfrak T}   \def\GTU {\mathfrak U}    \def\GTV {\mathfrak V}
\def\GTW {\mathfrak W}   \def\GTX {\mathfrak X}    \def\GTY {\mathfrak Y}
\def\GTZ {\mathfrak Z}   \def\GTQ {\mathfrak Q}

\font\Sym= msam10 
\def\SYM#1{\hbox{\Sym #1}}
\newcommand{\bdw}{\prt\Gw\xspace}
\date{}
\maketitle\medskip

\begin{abstract}
We study existence and stability of solutions of
 (E$_1$) $-\Delta   u +\frac{\mu}{|x|^{2 }}u+ g(u)= \gn\text{ in }\Omega,\;u=0\text{ on }\prt\Gw$,
 where $\Omega$ is a bounded, smooth domain of $\mathbb{ R}^N$, $N\geq  2$, containing the origin,
 $\mu\geq-\frac{(N-2)^2}{4}$ is a constant, $g$ is a nondecreasing function satisfying some integral
 growth assumption and $\gn$ is a Radon measure on  $\Gw$. We show that the situation differs according
 $\gn$ is diffuse or concentrated at the origin. When $g$ is a power we introduce a capacity framework to
 find necessary and sufficient condition for solvability.


\end{abstract}

\vspace{1mm}

\tableofcontents
\vspace{1mm}
\hspace{.05in}
\medskip

\setcounter{equation}{0}
\section{Introduction}

Let $\Omega\subset\mathbb{ R}^N$ ($N\geq 2$) be a bounded,  smooth domain containing the origin.
We define the elliptic operator with Leray-Hardy potential $\CL$ by
\bel{I}\mathcal{ L}_\mu:= -\Delta   +\frac{\mu}{|x|^{2 }},\ee
where $\gm$ is a real number satisfying
\begin{equation}\label{eq1.1a}
\gm\geq \gm_0:=-\frac{\left(N-2\right)^2}{4}.
\end{equation}
If $g:\R\to\R$ is a continuous nondecreasing function such that $g(0)\geq 0$, we are interested in the nonlinear Poisson equation
\begin{equation}\label{eq1.1}
\left\{
\begin{array}{lll}
 \mathcal{ L}_\mu  u+g(u)=\nu\quad & \rm{in}\ \,\Omega,\\[1.5mm]
 \phantom{  \mathcal{ L}_\mu    +g(u)}
u=0\quad & \rm{on}\  \partial \Omega,
\end{array}\right.
\end{equation}
where $\gn$ is a Radon measure in $\Gw$.

 When $\mu=0$, problem (\ref{eq1.1}) reduces  to
\begin{equation}\label{eq1.1-0}
\left\{
\begin{array}{lll}
-\Gd u +g(u)=\nu\quad & \rm{in}\ \,\Omega,\\[1.5mm]
 \phantom{  -\Gd    +g(u)}
u=0\quad & \rm{on}\  \partial \Omega,
\end{array}\right.
\end{equation}
which has been extensively studied by numerous authors in the
last 30 years. A fundamental contribution is due to Brezis
\cite{B12}, Benilan and Brezis \cite{BB11}, where $\nu$ is bounded and the function $g:\R\to\R$ is nondecreasing,
positive on $(0,+\infty)$ and satisfies the {\it  subcritical assumption} in dimension $N\geq 3$:
\bel{eq1.1-2}\int_1^{+\infty}(g(s)-g(-s))s^{-1-\frac{N}{N-2}}ds<+\infty.\ee
They obtained the existence, uniqueness and stability of weak solutions for the problem. When $N=2$,
 Vazquez \cite{Va1} introduced the exponential orders of growth of $g$ defined by
 \bel{I-x1}
\BA{lll}
\gb_+(g)=\inf\left\{b>0:\myint{1}{\infty}g\left(t\right)e^{-bt}dt<\infty\right\},\\[4mm]
\gb_-(g)=sup\left\{b<0:\myint{-\infty}{-1}g\left(t\right)e^{bt}dt>-\infty\right\}
\EA
 \ee
and proved that if  $\gn$ is  any bounded measure in $\Gw$ with Lebesgue decomposition
$$
\gn=\gn_r+\sum_{j\in\BBN}\ga_j\gd_{a_j},
$$
where $\gn_r$ is part of $\gn$ with no atom,  $a_j\in\Gw$ and $\ga_j\in\BBR$ satisfy
 \bel{I-x2}
\BA{lll}
\myfrac{4\gp}{\gb_-(g)}\leq \ga_j\leq \myfrac{4\gp}{\gb_+(g)},
\EA
 \ee
 then (\ref{eq1.1-0}) admits a (unique) weak solution.
Later on, Baras and Pierre \cite{BP2} studied (\ref{eq1.1-0})
when $g(u)=|u|^{p-1}u$ for $p>1$ and they discovered that if $p\geq\frac{N}{N-2}$
the problem is well posed if and only if $\nu$ is absolutely continuous
with respect to the Bessel capacity $c_{2,p'}$ with $p'=\frac{p}{p-1}$. \medskip


It is known that by the improved Hardy inequality \cite{BV1}  and Lax-Milgram theorem,   the non-homogeneous problem
 \begin{equation}\label{eq 1.1f}
 \displaystyle   \mathcal{L}_\mu u= f\qquad
   {\rm in}\quad  {\Omega},\qquad u=0\quad{\rm on}\quad \partial\Omega
\end{equation}
with $f\in L^2(\Omega)$,   has a unique solution  in $H^1_0(\Gw)$ if $\gm>\gm_0$, or in a weaker space $H(\Gw)$ if $\gm=\gm_0$ \cite{D}. A natural question is to find sharp condition of $f$ for the existence or nonexistence of (\ref{eq 1.1f}) and the difficulty comes from the fact that the Hardy term $| x|^{-2}u$ may not be locally integrable in $\Gw$.
An attempt done by Dupaigne in \cite{D} is to consider  problem (\ref{eq 1.1f}) when $\mu\in[\mu_0,0)$ and $N\geq 3$  in the sense of distributions
\begin{equation}\label{d1}
 \int_\Omega u\mathcal{L}_\mu \xi\, dx=\int_\Omega f\xi\, dx,\quad\forall\,\xi\in C^\infty_c(\Omega).
\end{equation}
The corresponding semi-linear problem is studied by \cite{BP} with this approach.\medskip

We adopt here a different point of view in using a different notion of weak solutions. It is known that the equation $\mathcal{ L}_\mu u=0$ in $\R^N\setminus \{0\}$ has two distinct radial solutions
 $$\Phi_\mu(x)=\left\{\arraycolsep=1pt
\begin{array}{lll}
 |x|^{\tau_-(\mu)}\quad
   &{\rm if}\quad \mu>\mu_0\\[1.5mm]
 \phantom{   }
|x|^{-\frac{N-2}{2}}\ln\left(\frac{1}{|x|}\right) \quad  &{\rm   if}\quad \mu=\mu_0
 \end{array}
 \right.\qquad \quad{\rm and}\qquad \Gamma_\mu(x)=|x|^{\tau_+(\mu)}$$
with
$$ \tau_-(\mu)=-\frac{N-2}2-\sqrt{\frac{(N-2)^2}4+\mu}\quad{\rm and}\quad   \tau_+(\mu)=-\frac{ N-2}2+\sqrt{\frac{(N-2) ^2}4+\mu}.$$
In the remaining of the paper and when there is no ambiguity,  we put $\tau_+=\tau_+(\mu)$, $\tau^0_+=\tau_+(\mu_0)$, $\tau_-=\tau_-(\mu)$ and $\tau^0_-=\tau_-(\mu_0)$. It is noticeable  that identity (\ref{d1}) cannot be used to express that $\Phi_\mu$ is a fundamental solution, i.e. $f=\gd_0$ since $\Phi_\mu$ is not locally integrable if $\gm\geq 2N$.  Recently, Chen, Quaas and Zhou found  in \cite{CQZ} that the function  $\Phi_\mu$ is the fundamental solution in the sense that it solves
\bel{L}\int_{\R^N}\Phi_\gm \mathcal{L}^*_\mu\xi \,d\gg_\mu(x)=c_{\mu}\xi(0),\quad \forall\xi\in C_c^{1,1}(\R^N),\ee
where
\begin{equation}\label{L1}
d{\gg_\mu}(x) =\Gamma_\mu(x) dx, \quad \ \mathcal{L}^*_\mu\xi=-\Delta\xi -2\frac{\tau_+ }{|x|^2}\langle x,\nabla\xi\rangle
\end{equation}
and
\begin{equation}\label{c-mu}
c_\gm=\left\{\BA {lll}2\sqrt{\gm-\gm_0}\abs{S^{N-1}}&\qquad\text{if }\;\gm>\gm_0,\\[2mm]
\left|S^{N-1}\right|&\qquad\text{if }\;\gm=\gm_0.
\EA\right.
\end{equation}

 With the power-absorption nonlinearity in $\Gw^*=\Omega\setminus\{0\}$, the precise behaviour near $0$ of any positive solution of
   \bel{I-1g}\mathcal{L}_\mu u+ u^p=0\qquad
   {\rm in}\ \, D'(\Gw^*)\ee
is given in \cite{GV}  when $p>1$. In this paper it appears a critical exponent
\begin{equation}\label{ce}
p_\mu^*=
1-\frac{2}{\tau_-},
\end{equation}
with the following properties: if $p\geq p_\mu^*$ any solution of (\ref{I-1g}) can be extended to be in $D'(\Gw)$. If $1<p<p_\mu^*$ any positive solution of (\ref{I-1g}) either satisfies
 \bel{I-1h}
 \lim_{x\to 0}|x|^{\frac2{p-1}}u(x)=\ell,
   \ee
 where $\ell=\ell_{N,p,\gm}>0$  or there exists $k\geq 0$ such that
    \bel{I-1i}
 \lim_{x\to 0}\myfrac{u(x)}{\Phi_\gm(x)}=k,
   \ee
and in that case  $u\in L^p_{loc}(\Gw;d\gg_\gm)$. In view of \cite{CQZ}, it implies that $u$ satisfies
\bel{I-1j}
\int_{\R^N}\left(u \mathcal{L}^*_\mu\xi +u^p \xi \right)d{\gg_\mu}(x)=c_{\mu}k \xi(0),\quad \forall\, \xi\in C_c^{1,1}(\R^N).
\ee
Furthermore, it is proved in \cite{GV} that when $\gm>\gm_0$ and $g:\R\to\R_+$ is a continuous nondecreasing function satisfying
\bel{I-2}
\myint{1}{\infty}\left(g(s)-g(-s)\right)s^{-1-p_\mu^*}ds<\infty,
\ee
then for any $k>0$ there exists a radial solution of
\bel{I-3}
\mathcal{ L}_\mu u+g(u)=0\qquad\text{in }\, D'(B^*_1)
 \ee
 satisfying (\ref{I-1i}), where $B^*_1:=B_1(0)\setminus\{0\}$.  When $\gm=\gm_0$ and $N\geq 3$ it is proved in \cite{GV} that if there exists $b>0$ such that
 \bel{I-2'}
\myint{0}{1} g\left(-b s^{-\frac{N-2}{N+2}}\ln s\right)ds<\infty,
\ee
then there a exists a radial solution of (\ref{I-3}) satisfying (\ref{I-1i}) with $\gg=\frac{(N+2)b}{2}$. In fact this condition is independent of $b>0$, by contrast the case $N=2$ and $\gm=0$ where the introduction of the exponential order of growth of $g$ is a necessity.
Moreover, it is easy to see that $u$ satisfies
 \bel{I-4}
\int_{\R^N}\left(u \mathcal{L}^*_\mu\xi + g(u)\xi\right)d{\gg_\mu}(x)=c_{\mu}\gg\xi(0),\quad \forall\xi\in C_c^{1,1}(\R^N).
\ee

In view of these results and identity (\ref{L}), we introduce a definition of weak solutions adapted to the operator
$\CL_\gm$ in a measure framework. Since $\Gg_\gm$ is singular at $0$ if $\gm<0$, there is  need of defining specific set of measures and {\it we denote by $\mathfrak{M}(\Omega^*; \Gamma_\mu)$, the set of Radon measures $\gn$  in $\Gw^*$ such that
 \bel{I-5}
 \myint{\Gw^*}{}\Gg_\gm d|\gn|:=\sup\left\{\myint{\Gw^*}{}\gz d|\gn|\, :\gz\in C_c(\Gw^*),\,0\leq\gz\leq\Gg_\gm\right\}<\infty.
\ee}
If $\gn\in \mathfrak{M}_+(\Omega^*)$, we define its natural extension, with the same notation since there is no ambiguity, as a measure in $\Omega$ by
 \bel{I-6}
\myint{\Gw}{}\gz d\gn=\sup\left\{\myint{\Gw^*}{}\eta  d\gn\, :\eta\in C_c(\Gw^*)\,,\;0\leq\eta\leq\gz\right\}\quad\text{for all }\gz \in C_c(\Gw)\,,\;\gz\geq 0,
\ee
a definition which is easily extended if $\gn=\gn_+-\gn_-\in \mathfrak{M}(\Omega^*)$. Since the idea is to use the weight $\Gg_\gm$ in the expression of the weak solution, the  expression $\Gg_\gm\gn$ has to be defined properly if $\gt_+<0$.
{\it We denote by $\mathfrak{M}(\Omega; \Gamma_\mu)$ the set of measures $\gn$ on $\Gw$ which coincide with the above natural extension of $\nu\lfloor_{\Gw^*}\in \mathfrak{M}_+(\Omega^*; \Gamma_\mu)$. } If $\gn\in \mathfrak{M}_+(\Omega; \Gamma_\mu)$ we define the measure $\Gg_\gm\gn$ in the following way
 \bel{I-6-1}
\myint{\Gw}{}\gz  d(\Gg_\gm\gn)=\sup\left\{\myint{\Gw^*}{}\eta  \Gg_\gm d\gn\, :\eta\in C_c(\Gw^*)\,,\;0\leq\eta\leq\gz\right\}\quad\text{for all }\gz \in C_c(\Gw)\,,\;\gz\geq 0.
\ee
If $\gn=\gn_+-\gn_-$, $\Gg_\gm\gn$ is defined acoordingly. Notice that the Dirac mass at $0$ does not belong to $\mathfrak{M}(\Omega; \Gamma_\mu)$ although it is a limit of $\{\gn_n\}\subset \mathfrak{M}(\Omega; \Gamma_\mu)$. We detote by $\overline{\mathfrak{M}}(\Omega; \Gamma_\mu)$ the set of measures which can be written under the form
\bel{I-7}
\nu=\gn\lfloor_{\Gw^*}+k\gd_0,
\ee
where $\gn\lfloor_{\Gw^*}\in  \mathfrak{M}(\Omega; \Gamma_\mu)$ and $k\in\R$.
Before stating our main
theorem we make precise the notion of weak solution used in this
article. We denote $\overline \Omega^*\!\!:=\overline\Gw\setminus\{0\}$, $\rho(x)={\rm dist}(x,\partial\Omega)$ and
\begin{equation}\label{test}
\BBX_\gm(\Gw)=\left\{\xi\in C_0(\overline\Gw)\cap C^1(\overline\Gw^*):|x|\CL_\gm^*\xi\in L^\infty(\Gw)\right\}.
\end{equation}
Clearly $C^{1,1}_0(\overline\Gw)\subset \BBX_\gm(\Gw)$.
\begin{definition}\label{weak definition}
We say that $u$ is a weak solution of (\ref{eq1.1}) with $\nu\in \overline{\mathfrak{M}}(\Omega; \Gamma_\mu)$ such that  $\nu=\gn\lfloor_{\Gw^*}+k\gd_0$ if $u\in
L^1(\Omega,\,|x|^{-1}d{\gg_\mu} )$, $g(u)\in L^1(\Omega, \gr d{\gg_\mu})$ and
\begin{equation}\label{weak sense}
\int_\Omega \left[u \mathcal{L}^*_\mu\xi+g(u)\xi\right]\,d{\gg_\mu}(x)=\myint{\Gw}{}\xi d(\Gg_\gm\nu)+k\xi(0)\quad \text{for all }\;\xi\in \BBX_\gm(\Gw),
\end{equation}
where $\mathcal{L}^*_\mu$ is given by (\ref{L}) and $c_{\mu}$ is defined in (\ref{c-mu}).

\end{definition}

A measure for which problem (\ref{eq1.1}) admits a solution is {\it a $g$-good measure}. In the regular case we prove the following\medskip

\nind{\bf Theorem A} {\it
Let $\mu\geq 0$  if  $N=2$, $\mu \geq   \mu_0$  if $ N\geq  3$ and $g:\R\to\R$ be a H\"{o}lder continuous nondecreasing function such that  $g(r)r\geq 0$ for any $ r\in\R$.
Then for any $\nu\in L^1(\Omega,d{\gg_\mu})$, problem (\ref{eq1.1})
has a unique  weak solution $u_\gn$ such that for some $c_1>0$,
$$\norm{u_\nu}_{L^1(\Omega,|x|^{-1} d{\gg_\mu})}\le c_1\norm{\nu}_{L^1(\Omega,d{\gg_\mu})}.$$
Furthermore, if $u_{\gn'}$ is the solution  of (\ref{eq1.1}) with right-hand side $\nu'\in L^1(\Omega,d{\gg_\mu})$, there holds
\bel{br1}
\myint{\Gw}{}\left[|u_{\gn}-u_{\gn'}|\CL^*_\gm\xi+|g(u_{\gn})- g(u_{\gn'})|\xi\right] d\gg_\gm(x)\leq \myint{\Gw}{}(\gn-\gn'){\rm sgn}(u-u')\xi d\gg_\gm(x)
\ee
and
\bel{br2}
\myint{\Gw}{}\left[(u_{\gn}-u_{\gn'})_+\CL^*_\gm\xi+(g(u_{\gn})- g(u_{\gn'}))_+\xi\right]d\gg_\gm(x)\leq \myint{\Gw}{}(\gn-\gn') {\rm sgn}_+(u-u')\xi d\gg_\gm(x)
\ee
for all $\xi\in \BBX_\gm(\Gw)$, $\xi\geq 0$.
}
\medskip

\begin{definition}\label{delta2} A continuous function  $g:\BBR\to\BBR$ such that $rg(r)\geq 0$ for all $r\in\BBR$ satisfies the weak $\Gd_2$-condition if there exists a positive nondecreasing function $t\in\R \mapsto K(t)$ such that
\bel{delta2}
|g(s+t)|\leq K(t)\left(|g(s)|+|g(t)|\right)\quad{\rm for\, all}\,\;(s,\,t)\in\BBR\ti\BBR\,\;{\rm s.t.}\,\;st\geq 0.
\ee
It satisfies the $\Gd_2$-condition if the above function $K$ is constant.
\end{definition}

Any power function or any exponential function satisfies the weak $\Gd_2$-condition.\medskip

\nind{\bf Theorem B} {\it
Let $\mu>0$ if  $N=2$ or $\mu >  \mu_0$ if $N\geq  3$ and $g:\BBR\to\BBR$ be a nondecreasing continuous function such that  $g(r)r\geq 0$ for any $ r\in\R$. If $g$ satisfies the weak $\Gd_2$-condition and
\begin{equation}\label{1.4-2}
  \int_1^\infty
(g(s)-g(-s))s^{-1-  \min\{p_\mu^*,\, p^*_0\}}ds<\infty,
\end{equation}
where $p_\mu^*$ is given by (\ref{ce}), then for any $\nu\in \overline{\mathfrak M}_+(\Gw;\Gg_\gm)$   problem (\ref{eq1.1}) admits a unique weak solution $u_\nu$.
}\medskip

Note that $\min\{p_\mu^*,\, p^*_0\}=p_\mu^*$ for $\mu>0$ and  $\min\{p_\mu^*,\, p^*_0\}=p_0^*$ if $\mu<0$. Furthermore, the
mapping: $\nu\mapsto u_\nu$ is increasing.
 In the case $N\geq 3$ and $\gm=\gm_0$ we have a more precise result.\medskip

\nind{\bf Theorem C} {\it
Assume that $N\geq  3$ and  $g:\BBR\to\BBR$ is a continuous nondecreasing function  such that  $g(r)r\geq 0$ for any $ r\in\R$ satisfying  the weak $\Gd_2$-condition and (\ref{eq1.1-2}).
Then for any $\nu={\nu\lfloor_{\Gw^*}}+c_{\mu}k\gd_0\in \overline{\mathfrak M}_+(\Gw;\Gg_\gm)$ problem (\ref{eq1.1}) admits a unique weak solution $u_\nu$.

Furthermore, if ${\nu\lfloor_{\Gw^*}}=0$, condition (\ref{eq1.1-2}) can be replaced by the following weaker one
\begin{equation}\label{1.4-3}
\myint{1}{\infty} \left(g(t)-g(-t)\right)\left(\ln t\right)^{\frac{N+2}{N-2}}t^{-\frac{2N}{N-2}}dt<\infty.
\end{equation}
}

Normally, the estimates on the Green kernel plays an essential role for approximating the solution of  elliptic problems with absorption and Radon measure. However, we  have banned the estimates on the Green kernel for Hardy operators due to luck them for $\mu\geq \mu_0$,  and our main idea is to separate the measure $\nu^*$  in $\mathfrak{M}(\Omega; \Gamma_\mu)$ and the Dirac mass at the origin,  and then to glue the solutions with above measures respectively.   This  requires   a very week new assumption: the weak $\Gd_2$-condition. \smallskip

In the previous result, it is noticeable that if $k=0$ (resp. $\nu\lfloor_{\Gw^*}=0$) only condition (\ref{eq1.1-2}) (resp. condition (\ref{1.4-3})) is needed. In the two cases the weak $\Gd_2$-condition is unnecessary.
In the power case where $g(u)=|u|^{p-1}u:=g_p(u)$,
\begin{equation}\label{1p.1}
\left\{
\begin{array}{lll}
 \mathcal{ L}_\mu  u+g_p(u)=\nu\ \ & \rm{in}\ \,\Omega,\\[2mm]
 \phantom{----\ \ \,\, }
u=0\ \ & \rm{on}\  \partial \Omega,
\end{array}\right.
\end{equation}
the following result follows from Theorem B and C.\medskip

\nind{\bf Corollary D} {\it Let $\gm\geq\gm_0$ if $N\geq 3$ and $\gm>0$ if $N= 2$. Any nonzero measure $\gn={\nu\lfloor_{\Gw^*}}+c_{\mu}k\gd_0\in \overline{\mathfrak M}_+(\Gw;\Gg_\gm)$ is $g_p$-good if  one of the following holds:  \smallskip

\nind (i)   $1<p<p^*_\gm$ in the case ${\nu\lfloor_{\Gw^*}}= 0$; \smallskip

\nind (ii)   $1<p<p^*_0$ in the case $k= 0$;\smallskip

\nind (iii)   $1<p<\min\left\{p^*_\gm,p^*_0\right\}$ in the case $k\neq 0$ and ${\nu\lfloor_{\Gw^*}}\neq 0$.

}\medskip

We remark that   $p^*_\gm$ is the sharp  exponent for existence of (\ref{1.4-3}) when ${\nu\lfloor_{\Gw^*}}= 0$,
while the critical exponent becomes $p^*_0$ when $k=0$ and $\nu$ has atom in $\Omega\setminus\{0\}$.\smallskip

The supercritical case of equation (\ref{1p.1}) corresponds to the fact that not all measures are $g_p$-good and the case where
$k\neq 0$ is already treated. \medskip


\nind{\bf Theorem E} {\it
Assume that $N\geq  3$. Then
 $\nu={\nu\lfloor_{\Gw^*}}\in \mathfrak{M}(\Omega; \Gamma_\mu)$ is $g_p$-good if and only if \smallskip
for any $\ge>0$, $\nu_\ge=\nu \chi_{_{B^c_\ge}}$ is absolutely continuous with respect to the $c_{2,p'}$-Bessel capacity.
}\medskip

Finally we characterize the compacts removable sets in $\Gw$.
 \medskip


\nind{\bf Theorem F} {\it
Assume that $N\geq  3$, $p>1$ and $K$ is a compact set of $\Gw$. Then any weak solution of
\begin{equation}\label{1p.2}
 \mathcal{ L}_\mu  u+g_p(u)= 0\quad  {\rm in}\ \, \Omega\setminus K
 \end{equation}
 can be extended a weak solution of the same equation in whole $\Gw$ if and only if \smallskip

\nind (i) $c_{2,p'}(K)=0$ if $0\notin K$;
 \smallskip

 \nind (ii) $p\geq p_{\gm^*}$ if $K=\{0\}$;\smallskip

 \nind (iii) $c_{2,p'}(K)=0$ if $\gm\geq 0$, $0\in K$ and $K\setminus\{0\}\neq\{\emptyset\}$;\smallskip

 \nind (iv) $c_{2,p'}(K)=0$ and $p\geq p_\gm^*$ if $\gm< 0$, $0\in K$ and $K\setminus\{0\}\neq\{\emptyset\}$.
}\medskip

 The case (i) is already proved in \cite[Theorem 1.2]{GV}. Notice also that if $A\neq\emptyset$ necessarily $c_{2,p'}(A)=0$ holds only if $p\geq p_0$. Therefore, if $\gm\geq 0$ there holds $p\geq p^*_0\geq p_\gm^*$, while if $\gm< 0$, then $p_0< p^*_\gm$.

\smallskip

  The rest of this  paper is organized as follows. In Section 2, we build the framework for weak solutions of (\ref{eq1.1}) involving $L^1$
  data. Section 3 is devoted to solve existence and uniqueness of weak solution  of (\ref{eq1.1}), where the absorption is subcritical and $\nu$ is a related Radon measure. Finally, we deal with the super critical case in Section 4 by characterized by Bessel Capacity.

\setcounter{equation}{0}
\section{$L^1$ data}


  Throughout this section we assume $N\geq 2$ and $\gm\geq\gm_0$ and in what follows, we denote by $c_i$ with $i\in\N$ a generic positive constant.  We first recall some classical comparison results for Hardy operator $\mathcal L_{\mu}$. The next lemma is proved in \cite[Lemma 2.1]{CQZ}, and  in  \cite[Lemma 2.1]{CC} if $h(s)=s^p$.

\begin{lemma}\label{lm cp}
Let $G$ be a bounded domain in $\R^N$ such that  $0\not\in \bar G$, $L: G\times [0,+\infty)\mapsto[0,+\infty)$ be a continuous function satisfying for any $x\in  G$,
$$h(x,s_1)\geq  h(x,s_2)\quad {\rm if}\quad s_1\geq  s_2 $$
and functions $u,\,v\in C^{1,1}(G)\cap C(\overline G)$ satisfy
$$\left\{
\BA {lll}\mathcal{L}_\mu u+ h(x,u)\geq  \mathcal{L}_\mu v+ h(x,v) \quad &{\rm in}\ \,  G,\\[1.5mm]
\phantom{\mathcal{L}_\mu + h(x,u)}   u\geq   v\quad &{\rm on}\ \partial G,
\EA\right.$$
then
$$u\geq  v\quad{\rm in}\quad  G.$$
\end{lemma}

As an immediate consequence we have

\begin{lemma}\label{cr hp}
Assume that $\Omega$ is a bounded $C^2$ domain containing $0$. If  $L$  is a continuous function as in Lemma \ref{lm cp} verifying furthermore
$L(x, 0)=0$ for all $x\in\Gw$, and $u\in C^{1,1}(\Omega^*)\cap C(\overline \Omega^*)$ satisfies
\begin{equation}\label{eq0 2.1}
 \arraycolsep=1pt\left\{
\begin{array}{lll} \phantom{ L}
 \displaystyle \mathcal{L}_\mu u  +L(x,u)= 0\quad\
   {\rm in}\  \,  \Omega^*,\\[1.5mm]\phantom{ L}
 \phantom{ L_\mu  +L(x,u)   }
 \displaystyle  u= 0\quad\ {\rm   on}\  \partial{\Omega},\\[1.5mm]
 \phantom{   }
  \displaystyle \lim_{x\to0}u(x)\Phi_\mu^{-1}(x)=0.
 \end{array}\right.
\end{equation}
Then $u=0$.
\end{lemma}

We recall that if $u\in L^1(\Omega, |x|^{-1}d{\gg_\mu})$ is a weak solution of
 \begin{equation}\label{homo}
 \arraycolsep=1pt\left\{
\begin{array}{lll}
  \mathcal{L}_\mu &u=f\quad\ & {\rm in}\ \,\Omega,\\[1mm]
&u=0\quad\ & {\rm on}\ \partial\Omega
\end{array}\right.
\end{equation}
in the sense of Definition \ref{weak definition}, it satisfies also
 \begin{equation}\label{weak-L1}
 \int_\Omega u \mathcal{L}^*_\mu(\xi)\, d{\gg_\mu}(x) =\int_\Omega f \xi\, d{\gg_\mu}(x)\quad{\rm for\,all} \ \xi\in \BBX_\gm(\Omega).\end{equation}
If $u$ is a weak solution of (\ref{homo}) there holds
 \begin{equation}\label{a}
\CL_\mu u=f\quad\ {\rm   in}\ \CD'(\Omega^*),
\end{equation}
and $v=\Gg_\gm^{-1}u$ verifies
 \begin{equation}\label{b}
\CL^*_\mu v=\Gg_\gm^{-1}f\quad\ {\rm   in}\ \CD'(\Omega^*).
\end{equation}

The following form of Kato's inequality, proved in \cite[Proposition 2.1]{CQZ}, plays an essential role in
the obtention a priori estimates and uniqueness of weak solution of (\ref{eq1.1}).

\begin{proposition}\label{pr 2.1}
If $f\in L^1(\Omega,\,\rho d{\gg_\mu})$, then there exists a  unique  weak
solution $u\in L^1(\Omega, |x|^{-1}d{\gg_\mu})$ of  (\ref{homo}).
Furthermore, for any $\xi\in  \BBX_\gm(\Omega)$, $\xi\geq 0$, we have
 \begin{equation}\label{sign}
\int_\Omega |u|  \mathcal{L}_\mu^*(\xi)\, d{\gg_\mu}(x) \le \int_\Omega
{{\rm sign}}(u)f  \xi\, d{\gg_\mu}(x)
\end{equation}
and
 \begin{equation}\label{sign+}
\int_\Omega u_+  \mathcal{L}_\mu^*(\xi)\,  d{\gg_\mu}(x) \le \int_\Omega
{{\rm sign}}_+(u)f \xi\, d{\gg_\mu}(x).
\end{equation}
\end{proposition}

\nind The proof is done if $\xi\in  C^{1,1}_0(\Omega)$, but it is valid if $\xi\in  \BBX_\gm(\Omega)$. The next result is proved in \cite[Lemma 2.3]{CZ}.

\begin{lemma}\label{lm 2.1-singular}
Assume that  $\mu>\mu_0$ and $f\in C^1(\Omega^*)$ verifies
\begin{equation}\label{2.1}
0\le f(x)\le c_2 |x|^{\tau-2}
\end{equation}
for some   $\tau>\tau_-$.
Let $u_f$ be the solution of
\begin{equation}\label{2.2}
 \arraycolsep=1pt\left\{
\begin{array}{lll} \phantom{ L----\ }
 \displaystyle \mathcal{L}_\mu u= f\quad\
   {\rm in}\ \, \Omega^*,\\[1.5mm]\phantom{ L}
 \phantom{ L_\mu- ---\ }
 \displaystyle  u= 0\quad\  {\rm   on}\  \partial{\Omega},\\[1.5mm]
 \phantom{}\!\!
  \displaystyle \lim_{x\to0}u(x)\Phi_\mu^{-1}(x)=0.
 \end{array}\right.
\end{equation}
Then there holds: \\[1mm]
(i) if $\tau_-<\tau<\tau_+$,
\begin{equation}\label{2.3}
0\le u_f(x)\le c_3|x|^{\tau}\ \quad    {\rm in}\quad \Omega^*;
\end{equation}
(ii) if $ \tau=\tau_+$,
\begin{equation}\label{2.3-1}
0\le u_f(x)\le c_4|x|^{\tau}(1+(-\ln |x|)_+)\ \quad    {\rm in}\quad \Omega^*;
\end{equation}
(iii) if $\tau>\tau_+$,
\begin{equation}\label{2.3-2}
0\le u_f(x)\le c_5|x|^{\tau_+} \ \quad    {\rm in}\quad \Omega^*.
\end{equation}

\end{lemma}

\nind{\it Proof of Theorem A}. Let  $\mathbb{H}^1_{\mu,0}(\Omega)$  be the closure of $C^\infty_0(\Omega)$ under the norm of
\bel{A1}\norm{u}_{\mathbb{H}^1_{\mu,0}(\Omega)}= \sqrt{\int_{\Omega} |\nabla u|^2 dx +\mu \int_{\Omega } \frac{ u^2}{|x|^2}  dx}.\ee
Then $\mathbb{H}^1_{\mu,0}(\Omega)$ is a Hilbert space with inner product
\bel{A2}\langle u, v \rangle_{\mathbb{H}^1_{\mu,0}(\Omega)}=\int_{\Omega} \langle \nabla u,\nabla v\rangle dx +\mu \int_{\Omega } \frac{ uv}{|x|^2}  dx\ee
and the embedding $\mathbb{H}^1_{\mu,0}(\Omega) \hookrightarrow L^p(\Omega)$ is continuous and compact
for $p\in[2, 2^*)$ with $2^*=\frac{2N}{N-2}$ when $N\geq 3$ and any $p\in [2],\infty$ if $N=2$. Furthermore, if $\eta\in C_c^1(\overline\Gw)$ has the value $1$ in a neighborhood of $0$, then  $\eta\Gg_\gm\in \mathbb{H}^1_{\mu,0}(\Omega)$. We put
$$G(v)=\myint{0}{v}g(s)ds,
$$
then $G$ is a convex nonnegative function. If $\gr\gn\in L^2(\Gw)$ we define the functional $J_\gn$  in the space $\mathbb{H}^1_{\mu,0}(\Omega)$ by
\bel{A3}J_\gn(v)=\left\{\BA {lll}\myfrac{1}{2}\norm{v}^2_{\mathbb{H}^1_{\mu,0}(\Omega)}+\myint{\Gw}{}G(v)dx-\myint{\Gw}{}\gn vdx&\quad {\rm if }\ G(v)\in L^1(\Gw,d\gg_\gm),\\[2mm]
\infty&\quad {\rm if }\ G(v)\notin L^1(\Gw,d\gg_\gm).
\EA\right.\ee
The functional $J$ is strictly convex, lower semicontinuous and coercive in $\mathbb{H}^1_{\mu,0}(\Omega)$, hence it admits a unique
minimum $u$ which satisfies
$$\langle u, v \rangle_{\mathbb{H}^1_{\mu,0}(\Omega)}+\myint{\Gw}{}g(u)v dx=\myint{\Gw}{}\gn v dx
\quad{\rm for\, all\;}v\in \mathbb{H}^1_{\mu,0}(\Omega).
$$
If $\xi\in C^{1,1}_0(\Gw)$ then $v=\xi\Gg_\gm\in \mathbb{H}^1_{\mu,0}(\Omega)$, then
\bel{A4}\BA {lll}\langle u, \xi\Gg_\gm \rangle_{\mathbb{H}^1_{\mu,0}(\Omega)}=
\myint{\Gw}{}\langle\nabla u,\nabla\xi\rangle d\gg_\gm(x)+\myint{\Gw}{}\left(\langle\nabla u,\nabla\Gg_\gm\rangle +\myfrac{\gm \Gamma_\mu }{|x|^2}\right)\xi dx
\EA\ee
and
$$\myint{\Gw}{}\langle\nabla u,\nabla\Gg_\gm\rangle\xi dx=-\myint{\Gw}{}\langle\nabla \xi ,\nabla\Gg_\gm\rangle udx-
\myint{\Gw}{}u\xi\Gd\Gg_\gm  dx,
$$
since $C^\infty_0(\Gw)$ is dense in $\mathbb{H}^1_{\mu,0}(\Omega)$. Furthermore, since $u\in L^p(\Gw)$ for any $p<2^*$,
$|x|^{-1}u\in L^1(\Gw,d\gg_\gm)$, hence $u\CL^*_\gm\xi\in L^1(\Gw,d\gg_\gm)$.
Therefore
\bel{A5}\myint{\Gw}{}\left(u\CL^*_\gm\xi+g(u)\xi\right)d\gg_\gm=\myint{\Gw}{}\gn\xi d\gg_\gm.
\ee
Next, if $\gn\in L^1(\Omega,\,\rho d{\gg_\mu})$ we consider a sequence $\{\gn_n\}\subset C^\infty_0(\Gw)$ converging to $\gn$ in
$L^1(\Omega,\,\rho d{\gg_\mu})$ and denote by $\{u_n\}$ the sequence of the corresponding minimizing problem in
$\mathbb{H}^1_{\mu,0}(\Omega)$. By Proposition \ref{pr 2.1} we have that, for any $\xi\in\BBX_\gm(\Gw)$,
\bel{A6}\myint{\Gw}{}\left(|u_n-u_m|\CL^*_\gm\xi+(g(u_n)-g(u_m)){\rm sgn}(u_n-u_m)\xi\right)d\gg_\gm
\leq \myint{\Gw}{}(\gn_n-\gn_m){\rm sgn}(u_n-u_m)\xi d\gg_\gm.
\ee
We denote by $\eta_0$ the solution of
\begin{equation}\label{A7}
\CL^*_\gm\eta=1\quad{\rm in}\ \ \Gw,\quad\ \eta=0\quad{\rm on}\ \ \prt\Gw.
\end{equation}
Its existence is proved in \cite[Lemma 2.2]{CQZ}, as well as the estimate $0\leq\eta_0\leq c_6\gr$ for some $c_6>0$. Since $g$ is monotone, we obtain from (\ref{A6})
\bel{A70}\myint{\Gw}{}\left(|u_n-u_m|+|g(u_n)-g(u_m)|\eta_0\right)d\gg_\gm
\leq \myint{\Gw}{}|\gn_n-\gn_m|\eta_0 d\gg_\gm.
\ee
Hence $\{u_n\}$ is a Cauchy sequence in $L^1(\Gw,d\gg_\gm)$.  Let $\eta_1$ solve the equation
\begin{equation}\label{A8}
\CL^*_\gm\eta=|x|^{-1}\quad\ {\rm in}\ \, \Gw^*,\qquad \eta=0\quad{\rm on}\ \, \prt\Gw.
\end{equation}
In the particular case $\Omega=B_1$,    function $ \eta_1(x)=\frac{1-|x|}{(N-1+2\gt_+(\gm))}$ verifies
$$ \arraycolsep=1pt\left\{ \BA {lll}\CL^*_\gm\eta_1= |x|^{-1} \quad\ &{\rm in}\ \,B^*_1,\\[1.5mm]
\phantom{\CL^*_\gm}
\eta_1=0\quad&{\rm on}\ \,\prt B_1
\EA\right.$$
(we can always assume that $\Gw\subset B_1$). As in the proof of \cite[Lemma 2.2]{CQZ}, for any $x_0\in \Gw$ there exists
$r_0>0$ such that $B_{r_0}(x_0)\subset\Gw$ and for $t>0$ small enough $w_{t,x_0}(x)=t(r_0^2-|x-x_0|^2)$ is a subsolution of
(\ref{A7}), hence of (\ref{A8}). Therefore $\eta_1$ exists. Using again the density of $C^\infty_0(\Gw)$ in $\mathbb{H}^1_{\mu,0}(\Omega)$ and integrating on $\Gw\setminus B_\ge$ and letting $\ge\to 0$, we obtain as a variant of (\ref{A70})
\bel{A81}\myint{\Gw}{}\left(\myfrac{|u_n-u_m|}{|x|}+|g(u_n)-g(u_m)|\eta_1\right)d\gg_\gm(x)
\leq \myint{\Gw}{}|\gn_n-\gn_m|\eta_1 d\gg_\gm.
\ee
Hence $\{u_n\}$ is a Cauchy sequence in $L^1(\Omega,\,|x|^{-1} d{\gg_\mu})$ with limit $u$ and $\{g(u_n)\}$ is a Cauchy sequence in $L^1(\Omega,\,\rho d{\gg_\mu})$ with limit $g(u)$.
Then (\ref{A5}) holds. As for (\ref{br1}) it is a consequence of (\ref{A6}) and (\ref{br2}) is proved similarly.\qeda
\smallskip

\setcounter{equation}{0}
\section{The subcritical case}

In this section as well as in the next one we always assume that $N\geq 3$ and $\gm\geq\gm_0$, or  $N= 2$ and $\gm>0$, since the case $N=2$, $\gm=0$, which necessitates specific tools, has already been completely treated in \cite{Va1}.

We recall that the set $\frak M(\Gw^*;\Gg_\gm)$ of Radon measures is defined in introduction as the set of measures in  $\Gw^*$ satisfying (\ref{I-5}), and any positive measure $\gn\in\frak M(\Gw^*;\Gg_\gm)$ is naturaly extended by formula (\ref{I-6}) as a positive measure in $\Gw$. The space $\overline{\frak M}(\Gw;\Gg_\gm)$ is the space of measures $\gn$ on $C_c(\Gw)$ such that
\bel{DM1}\displaystyle
\gn=\gn\lfloor_{\Gw^*}+k\gd_0,
\ee
where $\gn\lfloor_{\Gw^*}\in \frak M(\Gw^*;\Gg_\gm)$.
\subsection{The linear equation}
\blemma{line} If $\gn\in \overline{\frak M}(\Gw;\Gg_\gm)$, then there exists a unique weak solution $u\in L^1(\Gw,|x|^{-1}d\gg_\mu)$ to
\begin{equation}\label{line1}
 \arraycolsep=1pt\left\{
\BA {lll}
\CL_\gm u=\gn\quad & \rm{in}\;\,\Omega,\\[1.5mm]
\phantom{\CL_\gm }u=0\quad & \rm{on}\;\,\prt\Omega.
\EA\right.
\end{equation}
This solution is denoted by $\BBG_{\gm}[\gn]$, and this defines the Green operator of $\CL_\gm$ in $\Gw$ with homogeneous Dirichlet conditions.
\es
\Proof By linearity and using the result of \cite{CQZ} on fundamental solution, we can assume that $k=0$ and $\gn\geq 0$. Let $\{\gn_n\}\subset L^1(\Gw,\gr d\gg_\mu)$ be a sequence
 such that $\gn_n\geq 0$ and
$$\myint{\Gw}{}\xi\Gg_\gm\gn_n dx\to \myint{\Gw}{}\xi d(\Gg_\gm\gn)\quad \rm {for\, all}\;\,\xi\in \BBX_\gm(\Gw),
$$
and by Proposition \ref{pr 2.1}, we may let $ u_n$ be the unique, nonnegative weak solution of
\begin{equation}\label{line2} \arraycolsep=1pt\left\{
\BA {lll}
\CL_\gm u_n=\gn_n\quad & \rm{in}\;\,\Omega,\\[1.5mm]
\phantom{\CL_\gm }u_n=0\quad & \rm{on}\;\,\prt\Omega
\EA\right.
\end{equation}
with $n\in \N$.
  There holds
\begin{equation}\label{line2'}\BA {lll}
\myint{\Gw}{}u_n\CL^*_\gm\xi d\gg_\gm(x)=\myint{\Gw}{}\xi \gn_n\Gg_\gm dx\ \quad \rm{for\,all}\;\,\xi\in\BBX_\gm(\Omega).
\EA
\end{equation}
Then $u_n\geq 0$ and using the function $\eta_1$ defined in the proof of Theorem A for test function, we have
\begin{equation}\label{line3}\BA {lll}
c\myint{\Gw}{}\myfrac{u_n}{|x|}d\gg_\gm=\myint{\Gw}{}\eta_1\Gg_\gm\gn_n dx\le c\norm{  \nu }_{\mathfrak M(\Omega,\Gamma_\mu )} ,
\EA
\end{equation}
which implies that $\{u_n\}$ is bounded in $L^1(\Omega, \frac1{|x|}d\gg_\gm(x))$.

For any $\epsilon>0$ sufficiently small, set the test function $\xi$ in $\left\{\zeta\in \BBX_\gm(\Omega):\zeta=0\ \, {\rm in}\ B_\epsilon\right\}$, then we have that
\begin{equation}\label{line2'-}\BA {lll}
\myint{\Gw\setminus B_\epsilon(0)}{}u_n\CL^*_\gm\xi d\gg_\gm(x)=\myint{\Gw\setminus B_\epsilon(0)}{}\xi \gn_n\Gg_\gm dx\ \quad \rm{for\,all}\;\,\xi\in\BBX_\gm(\Omega).
\EA
\end{equation}
Therefore,  for any open sets $O,O'$ verifying  $\bar O \subset  O'\subset \bar O'\subset \Omega\setminus B_\epsilon(0)$, there exists $c>0$ independent of $n$ such that
$$\norm{u_n}_{L^1(O')}\le c\norm{  \nu }_{\mathfrak M(\Omega,\Gamma_\mu )} . $$
 Note that in $\Omega\setminus B_\epsilon$, the operator $\CL^*_\gm$ is uniformly elliptic  and the measure $d\gg_\gm$
 is equivalent to the Hausdorff measure $dx$,
then  \cite[Corollary 2.8]{V}  could be applied to obtain that for some $c>0$ independent of $n$ but dependent of $O'$,
 \begin{eqnarray*}
\norm{u_n}_{W^{1,q}(O)} &\le & c\norm{u_n}_{L^1(O')}+\norm{\tilde \nu_{n}}_{L^1(\Omega,d\gg_\gm)}  \\
    &\le &c \norm{  \nu }_{\mathfrak M (\Omega,\Gamma_\mu )}.
 \end{eqnarray*}
That is, $\{u_n\}$ is uniformly bounded in $W^{1,q}_{loc}(\Omega\setminus\{0\})$.

As a consequence, by the arbitrary of $\epsilon$, there exist a subsequence, still denoting  $\{u_{n}\}_n$ and $u$  such that
$$u_{n}\to u\quad{\rm a.e.\ in}\ \Omega.$$
Meanwhile, we deduce from Fatou's lemma,
\begin{equation}\label{line4}\BA {lll}
\myint{\Gw}{}\myfrac{u}{|x|}d\gg_\mu\leq  c\myint{\Gw}{}\eta_1 \Gg_\gm d\gn.
\EA
\end{equation}

We next claim  that  $u_n\to u$ in $L^1(\Gw,|x|^{-1}d\gg_\mu)$. Let $\gw\subset\Gw$ be a Borel set and $\psi_\gw$ be the solution of
\begin{equation}\label{line5}
\arraycolsep=1pt\left\{\BA {lll}
\CL^*_\gm \psi_\gw=|x|^{-1}\chi_{{\gw}}\quad & {\rm in}\ \,\Omega,\\[1.5mm]
\phantom{\CL^*_\gm }\psi_\gw=0\quad & {\rm on}\ \prt\Omega.
\EA\right.
\end{equation}
Then $\psi_\gw\leq \eta_1$, thus it is uniformly bounded. Assuming that $\Gw\subset B_1$, clearly $\psi_\gw$ is bounded from above by the solution $\Psi_\gw$ of
\begin{equation}\label{line6}
 \arraycolsep=1pt\left\{
\BA {lll}
\CL^*_\gm \Psi_\gw=|x|^{-1}\chi_{{\gw}}\quad & {\rm in}\;\,B_1,\\[1.5mm]
\phantom{\CL^*_\gm }\Psi_\gw=0\quad & {\rm on}\;\,\prt B_1
\EA\right.
\end{equation}
and by standard rearrangement, $\sup_{B_1}\Psi_\gw\leq \sup_{B_1}\Psi^r_\gw$, where $\Psi^r_\gw$ solves
\begin{equation}\label{line7}
 \arraycolsep=1pt\left\{\BA {lll}
\CL^*_\gm \Psi^r_\gw=|x|^{-1}B_{\ge(|\gw|)}\quad & {\rm in}\;\,B_1,\\[1.5mm]
\phantom{\CL^*_\gm } \Psi^r_\gw=0\quad & {\rm on}\ \prt B_1,
\EA\right.
\end{equation}
where $\ge(|\gw|)=\left(\frac{|\gw|}{|B_1}\right)^{\frac{1}{N}}$. Then $\Psi^r_\gw$ is radially decreasing and $\lim_{|\gw|\to 0}\Psi^r_\gw=0$, uniformly on $B_1$. This implies
\begin{equation}\label{line8}\BA {lll}\displaystyle
\lim_{|\gw|\to 0}\psi_\gw(x)=0\quad  {\rm uniformly\ in}\ \,B_1.
\EA
\end{equation}
Using (\ref{line2'}) with $\xi=\psi_\gw$,
$$\displaystyle\myint{\gw}{}\myfrac{u_n}{|x|}d\gg_\gm(x)=\myint{\gw}{}\gn_n\Gg_\gm\psi_\gw dx\leq \sup_{\Gw}\psi_\gw\myint{\gw}{}\gn_n\Gg_\gm dx\to 0\ \ \rm{as}\;\,|\gw|\to 0.
$$
Therefore $\{u_n\}$ is uniformly integrable for the measure $|x|^{-1}d\gg_\gm$. Letting $n\to\infty$ in (\ref{line2'}) implies the claim.\qeda
\medskip

\subsection{Dirac masses}

We assume that $g:\BBR\to\BBR$ is a continuous nondecreasing function such that $rg(r)\geq 0$ for all $r\in\BBR$.
The next lemma dealing with problem
\begin{equation}\label{sc1}
\arraycolsep=1pt\left\{
\begin{array}{lll}
 \mathcal{ L}_\mu  u+g(u)=k\gd_0\quad & \rm{in}\ \,\Omega,\\[1.5mm]
 \phantom{  \mathcal{ L}_\mu    +g(u)}
u=0\quad & \rm{on}\  \partial \Omega
\end{array}\right.
\end{equation}
 is an extension of \cite[Theorem 3.1, Theorem 3.2]{GV}. Actually it was quoted in this article as Remark 3.1 and Remark 3.2 and we give here their proof. Notice also that when $N\geq 3$ and $\gm=\gm_0$ we give a more complete result that \cite[Theorem 3.2]{GV}.

\blemma{dirac} Let $k\in\BBR$ and $g:\BBR\to\BBR$ be a continuous nondecreasing function such that $rg(r)\geq 0$ for all $r\in\BBR$. Then problem (\ref{sc1}) admits a unique solution $u:=u_{k\gd_0}$ if   one of the following conditions is satisfied:\smallskip

\nind (i) $N\geq 2$, $\gm > \gm_0$ and $g$ satisfies (\ref{I-2});\smallskip

\nind (ii)  $N\geq 3$, $\gm=\gm_0$ and $g$ satisfies (\ref{1.4-3}).
\es
\Proof Without loss of generality we assume $B_R\subset \Gw\subset B_1$  for some $R\in (0,1)$.\\
\nind (i) {\it The case  $\gm > \gm_0$}. It follows from
\cite[Theorem 3.1]{GV} that for any $k\in\BBR$ there exists a radial function  $v_{k,1}$ (resp. $v_{k,R}$)  defined in $B_1^*$ (resp. $B_R^*$) satisfying
\begin{equation}\label{sc2}
\CL_\gm v +g(v)=0\qquad{\rm in}\;B_1^*\quad {\rm(resp. \; in }\;B_R^*{\rm )},
\end{equation}
vanishing respectively on $\prt B_1$ and $\prt B_R$ and satisfying
\begin{equation}\label{sc3}\displaystyle
\lim_{x\to 0}\myfrac{v_{k,1}(x)}{\Phi_\gm(x)}=\lim_{x\to 0}\myfrac{v_{k,R}(x)}{\Phi_\gm(x)}=\myfrac{k}{c_\gm}.
\end{equation}
Furthermore $g(v_{k,1})\in L^1(B_1,d\gg_\gm)$ (resp. $g(v_{k,R})\in L^1(B_R,d\gg_\gm)$). Assume that $k>0$, then $0\leq v_{k,R}\leq v_{k,1}$ in $B_R^*$ and the extension of $\tilde v_{k,R}$ by $0$ in $\Gw^*$ is a subsolution of (\ref{sc2}) in $\Gw^*$ and it is still smaller than $v_{k,1}$ in $\Gw^*$. By the well known method on super and subsolutions (see e.g. \cite[Theorem 1.4.6]{Vebook}), there exists a function $u$ in $\Gw^*$ satisfying $\tilde v_{k,R}\leq u\leq v_{k,1}$ in $\Gw^*$ and
\begin{equation}\label{sc4}
\arraycolsep=1pt\left\{\BA {lll}
\CL_\gm u +g(u)=0& \quad{\rm in}\,\ \Gw^*,\\[1mm]
\phantom{\CL_\gm  +g(u)}
u=0&  \quad{\rm on} \ \prt\Gw,\\[1mm]
\phantom{,}
\displaystyle
\lim_{x\to 0}\myfrac{u(x)}{\Phi_\gm(x)}=\myfrac{k}{c_\gm}.
\EA\right.
\end{equation}
By standard methods in the study of isolated singularities (see e.g. \cite {GV}, \cite{Vebook0},  and \cite {C} and \cite {CD} for various extensions)
\begin{equation}\label{sc5}\displaystyle
\lim_{x\to 0}|x|^{1-\gt_-}\nabla u(x)=\gt_-\myfrac{k}{c_\gm}\myfrac{x}{|x|}.
\end{equation}
For any $\ge>0$ and $\xi\in\BBX_\gm(\Gw)$,
$$\BA {lll}0=\myint{\Gw\setminus B_\ge}{}(\CL_\gm u +g(u))\Gg_\gm\xi dx\\[4mm]\phantom{0}
=\myint{\Gw\setminus B_\ge}{}u\CL^*_\gm \xi d\gg_\gm(x)+(\gt_--\gt_+)\myfrac{k}{c_\gm}|S^{N-1}|\xi(0)(1+o(1)).
\EA$$
Using (\ref{c-mu}), we obtain
\begin{equation}\label{sc6}\displaystyle
\myint{\Gw}{}u\CL^*_\gm \xi d\gg_\gm(x) =k\xi(0).
\end{equation}

\nind (ii){\it The case  $\gm = \gm_0$}. In \cite[Theorem 3.2]{GV} it is proved that if for some $b>0$ there holds
\begin{equation}\label{sc6-1}\displaystyle
I:=\myint{1}{\infty}g\left(bt^{\frac{N-2}{N+2}}\ln t\right) t^{-2} dt<\infty,
\end{equation}
then there exists a solution of (\ref{I-3}) satisfying (\ref{I-1i}) with $\gg=\frac{(N+2)b}{2}$.  Actually we claim that {\it  the finiteness of this integral is independent of the value of $b$}. To see that, set $s=t^{\frac{N-2}{N+2}}$, then
$$I=\frac{N+2}{N-2}\myint{1}{\infty}g\left(\gb s\ln s\right) s^{-\frac{2N}{N-2}} ds
$$
with $\gb=\frac{N+2}{N-2}b$. Set $\gt=\gb s\ln s$, then
$$\ln s\left(1+\frac{\ln\ln s}{\ln s}+\frac{\ln\gb}{\ln s}\right)\Longrightarrow \ln s=\ln\gt(1+o(1))\quad{\rm as}\;\,s\to\infty.$$
We infer that for $\ge>0$ there exists $s_\ge>2$ and $\gt_\ge=s_\ge\ln s_\ge$ such that
\begin{equation}\label{sc6-2}
(1-\ge)\gb^{\frac{N+2}{N-2}}\leq
\myfrac{\myint{s_\ge}{\infty}g\left(\gb s\ln s\right) s^{-\frac{2N}{N-2}} ds}{\myint{\gt_\ge}{\infty}g\left(\gt\right)
\left(\ln\gt\right)^{\frac{N+2}{N-2}}\gt^{-\frac{2N}{N-2}} d\gt}\leq (1+\ge)\gb^{\frac{N+2}{N-2}},
\end{equation}
which implies the claim. Next we prove as in case (i) the existence of $v_{k,1}$ (resp. $v_{k,R}$)  defined in $B_1^*$ (resp. $B_R^*$) satisfying
\begin{equation}\label{sc7}
\CL_{\gm_0} v +g(v)=0\qquad{\rm in}\;B_1^*\quad {\rm(resp. \; in }\;B_R^*{\rm )},
\end{equation}
vanishing respectively on $\prt B_1$ and $\prt B_R$ and satisfying
\begin{equation}\label{sc8}\displaystyle
\lim_{x\to 0}\myfrac{v_{k,1}(x)}{\Phi_\gm(x)}=\lim_{x\to 0}\myfrac{v_{k,R}(x)}{\Phi_\gm(x)}=\myfrac{k}{c_{\gm_0}}.
\end{equation}
We end the proof as above.\qeda

\medskip

\nind\Remark It is important to notice that conditions (\ref{I-2}) and (\ref{1.4-3}) (or equivalently (\ref{I-2'})) are also necessary for the existence of radial solutions in a ball, hence their are also necessary for the existence of non radial solutions of the Dirichlet problem (\ref{sc1}).
\subsection{Measures in $\Gw^*$}
We consider now the problem
\begin{equation}\label{sc9}
 \arraycolsep=1pt\left\{
\begin{array}{lll}
 \mathcal{ L}_\mu  u+g(u)=\gn\quad & \rm{in}\ \,\Omega,\\[1.5mm]
 \phantom{  \mathcal{ L}_\mu    +g(u)}
u=0\quad & \rm{on}\  \partial \Omega,
\end{array}\right.
\end{equation}
where  $\gn\in\mathfrak M(\Gw^*;\Gg_\gm)$.

\begin{lemma}\label{noatom} Let $\gm\geq\gm_0$. Assume that $g$ satisfies (\ref{eq1.1-2}) if $N\geq 3$ or the $\gb_\pm(g)$ defined by
(\ref{I-x1}) satisfy $\gb_-(g)<0<\gb_+(g)$ if $N=2$, and let $\gn\in\mathfrak M(\Gw^*;\Gg_\gm)$. If $N=2$, we assume that
$\gn$ can be decomposed as $\gn=\gn_r+\sum_j\ga_j\gd_{a_j}$ where $\gn_r$ has no atom,  the $\ga_j$ satisfy (\ref{I-x2}) and $\{a_j\}\subset\Gw^*$. Then
problem (\ref{sc9}) admits a unique weak solution.
\end{lemma}
\Proof We assume first that $\gn\geq 0$ and let $r_0=\dist (x,\prt\Gw)$. For $0<\gs<r_0$, we set $\Gw^\gs=\Gw\setminus\{\overline B_\gs\}$ and $\gn_\gs=\gn\chi_{_{\Gw^\gs}}$ and for $0<\ge<\gs$ we consider the following problem in $\Gw^\ge$
\begin{equation}\label{sc10}
 \arraycolsep=1pt\left\{
\begin{array}{lll}
 \mathcal{ L}_\mu  u+g(u)=\gn_\gs\quad & \rm{in}\ \,\Omega^\ge,\\[1mm]
 \phantom{  \mathcal{ L}_\mu    +g(u)}
u=0\quad & \rm{on}\  \partial \Omega,\\[1mm]
 \phantom{  \mathcal{ L}_\mu    +g(u)}
 u=0& \rm{on}\  \partial B_\ge.
\end{array}\right.
\end{equation}
Since $0\notin\Gw^\ge$   problem (\ref{sc10}) admits a unique solution $u_{\gn_\gs,\ge}$ which is smaller than $\BBG_{\gm}[\gn]$ and satisfies
$$0\leq u_{\gn_\gs,\ge}\leq u_{\gn_{\gs'},\ge'}\quad  \rm{in}\;\,\Omega^\ge\;\,  \rm{for\, all}\;\, 0<\ge'\leq\ge\,\;\rm{and}\;\,0<\gs'\leq\gs.
$$
For any $\xi\in \C^{1,1}_c(\Gw^*)$ and $\ge$ small enough so that supp$\,(\xi)\subset \Gw^\ge$, there holds
$$\myint{\Gw}{}\left(u_{\gn_\gs,\ge}\CL_\gm^*\xi+g(u_{\gn_\gs,\ge})\xi\right) d\gg_\mu =\myint{\Gw}{}\xi \Gg_\mu d\gn_\gs.
$$
There exists $\displaystyle u_{\gn_\gs}=\lim_{\ge\to 0}u_{\gn_\gs,\ge}$ and it satisfies the identity
\begin{equation}\label{sc11}
\myint{\Gw}{}\left(u_{\gn_\gs}\CL_\gm^*\xi+g(u_{\gn_\gs})\xi\right) d\gg_\mu =\myint{\Gw}{}\xi \Gg_\mu d\gn_\gs\quad\rm{for\, all}\;\,\xi\in \C^{1,1}_c(\Gw^*).
\end{equation}
Using the maximum principle and \rlemma{line}, there holds
\begin{equation}\label{sc12}
0\leq u_{\gn_\gs}\leq \BBG_{\gm}[\gn_\gs]\leq \BBG_{\gm}[\gn].
\end{equation}
Since $\gn_\gs$ vanishes in $B_\gs$, $\BBG_{\gm}[\gn_\gs](x)\leq c\Phi_\gm(x)$ in a neighborhood of $0$, and $u_{\gn_\gs}$ is also bounded by $c\Phi_\gm$ in this neighborhood. This implies that
$\Phi^{-1}_\gm(x)u_{\gn_\gs}(x)\to c'$ as $x\to 0$ for some $c'\geq 0$. Next let $\xi\in \C^{1,1}_c(\Gw)$,
$$\ell_n(r)=\left\{\BA {lll}
2^{-1}\left(1+\cos\left(\frac{2\gp |x|}{\gs}\right)\right)\quad &{\rm if}\;\ \,|x|\leq\frac{\gs}{2},\\[1.5mm]
0&{\rm if}\;\,\ |x|>\frac{\gs}{2}
\EA\right.$$
and $\xi_n=\xi\ell_n$. Then
\begin{equation}\label{sc13}
\myint{\Gw}{}\left(u_{\gn_\gs}\CL_\gm^*\xi_n+g(u_{\gn_\gs})\xi_n\right) d\gg_\mu =\myint{\Gw}{}\xi_n \Gg_\mu d\gn_\gs.
\end{equation}
When $n\to\infty$,
$$\myint{\Gw}{}\xi_n \Gg_\mu d\gn_\gs\to \myint{\Gw}{}\xi \Gg_\mu d\gn_\gs
$$
and
$$\myint{\Gw}{}g(u_{\gs})\xi_nd\gg_\mu \to\myint{\Gw}{}g(u_{\gs})\xi d\gg_\mu.
$$
Now for the first inegral term in (\ref{sc13}), we have
$$\myint{\Gw}{}u_{\gn_\gs}\CL_\gm^*\xi_n d\gg_\mu=\myint{\Gw}{}\ell_nu_{\gs}\CL_\gm^*\xi d\gg_\mu+I_n+II_n+III_n,
$$
where
$$I_n=-\myint{B_{\frac\gs 2}}{}u_{\gs}\xi\Gd\ell_n d\gg_\mu,
$$
$$II_n=-2\myint{B_{\frac\gs 2}}{}u_{\gs}\langle\nabla\xi,\nabla\ell_n\rangle d\gg_\mu
$$
and
$$III_n=-\gt_+\myint{B_{\frac\gs 2}}{}u_{\gs}\langle \frac{x}{|x|^2},\nabla\ell_n\rangle d\gg_\mu.
$$
Using the fact that $\xi(x)\to \xi(0)$ and $\nabla\xi(x)\to \nabla\xi(0)$ we easily infer that $I_n$, $II_n$ and $III_n$ to 0 when
$n\to\infty$, the most complicated case being the case when $\gm=\gm_0$, which is the justification of introducing the explicit cut-off function $\ell_n$. Therefore (\ref{sc11}) is still valid if it is assumed that $\xi\in \C^{1,1}_c(\Gw)$. This means that $u_{\gn_\gs}$ is a weak solution of
\begin{equation}\label{sc14}
\left\{
\begin{array}{lll}
 \mathcal{ L}_\mu  u+g(u)=\gn_\gs\quad & \rm{in}\ \,\Omega,\\[1mm]
 \phantom{  \mathcal{ L}_\mu    +g(u)}
u=0\quad & \rm{on}\  \partial \Omega.
\end{array}\right.
\end{equation}
Furthermore $u_{\gn_\gs}$ is unique and $u_{\gn_\gs}$ is a decreasing function of $\gs$ with limit $u$ when $\gs\to 0$. Taking $\eta_1$ as test function, we have
$$\myint{\Gw}{}\left(c|x|^{-1}u_{\gn_\gs} +\eta_1g(u_{\gn_\gs})\right)d\gg_\mu=\myint{\Gw}{}\eta_1d\left(\gg_\mu\gn_\gs\right)\leq \myint{\Gw}{}\eta_1d\left(\gg_\mu\gn\right).
$$
Using monotone convergence theorem we infer that $u_{\gn_\gs}\to u$ in $L^1(\Gw,|x|^{-1}d\gg_\gm)$ and $g(u_{\gn_\gs})\to g(u_\gn)$ in $L^1(\Gw,d\gg_\gm)$. Hence $u=u_\gn$ is the weak solution of (\ref{sc9}).\smallskip

\nind Next we consider a signed measure $\gn=\gn_+-\gn_-$. We denote by $u_{\gn^\gs_+,\ge}$, $u_{-\gn^\gs_-,\ge}$ and $u_{\gn^\gs,\ge}$ the solutions of (\ref{sc10}) in $\Gw^\ge$ corresponding to $\gn^\gs_+$, $-\gn^\gs_-$ and $\gn^\gs,\ge$ respectively. Then
\begin{equation}\label{sc15}
u_{-\gn^\gs_-,\ge}\leq u_{\gn^\gs,\ge}\leq  u_{\gn^\gs_+,\ge}.
\end{equation}
The correspondence $\ge\mapsto u_{\gn^\gs_+,\ge}$ and $\ge\mapsto u_{-\gn^\gs_-,\ge}$ are respectively increasing and decreasing. Furthermore $u_{\gn^\gs,\ge}$ is locally bounded, hence by local compactness and up to a subsequence $u_{\gn^\gs,\ge}$ converges a.e. in $B_\ge$ to some function $u_{\gn^\gs}$. Since $u_{-\gn^\gs_-,\ge}\to u_{-\gn^\gs_-}$ and $u_{\gn^\gs_+,\ge}\to u_{\gn^\gs_+}$ in $L^1(\Gw,|x|^{-1}d\gg_\gm)$, it follows by Vitali's theorem that $u_{\gn^\gs,\ge}\to u_{\gn^\gs}$ in $L^1(\Gw,|x|^{-1}d\gg_\gm)$. Similarly, using the monotonicity of $g$, $g(u_{\gn^\gs,\ge})\to g(u_{\gn^\gs})$ in $L^1(\Gw,d\gg_\gm)$. By local compactness, $u_{\gn^\gs}\to u$ a.e. in $\Gw$. Using the same argument of uniform integrability, we have that $u_{\gn^\gs}\to u$ in $L^1(\Gw,|x|^{-1}d\gg_\gm)$ and $g(u_{\gn^\gs})\to g(u)$ in $L^1(\Gw,d\gg_\gm)$ when $\gs\to 0$ and $u$ satisfies
\begin{equation}\label{sc16}
\myint{\Gw}{}\left(u\CL_\gm^*\xi+g(u)\xi\right) d\gg_\mu =\myint{\Gw}{}\xi d(d\gg_\mu\gn)\quad{\rm for\ any} \ \xi\in \C^{1,1}_c(\Gw^*).
\end{equation}
 Finally the singularity at $0$ is removable by the same argument as above which implies that $u$ solves (\ref{sc16}) and thus $u=u_\gn$ is the weak solution of (\ref{sc9}).\qeda
\subsection{Proof of Theorem B}

The idea is to glue altogether two solutions one with the Dirac mass and the other with the measure in $\Gw^*$, this is the reason why the weak $\Gd_2$ condition is introduced.
\blemma{sigma} Let $\gn=\gn\lfloor_{\Gw^*}+k\gd_0\in\overline{\frak M}_+(\Gw;\Gg_\gm)$ and $\gs>0$. We assume that $\gn\lfloor_{\Gw^*}(\overline B_\gs)=0$. Then there exists a unique weak solution to (\ref{eq1.1}).
\es
\Proof Set $\gn_\gs=\gn\lfloor_{\Gw^*}$. It has support in $\Gw_\gs=\Gw\setminus\overline B_\gs$. For $0<\ge<\gs$ we consider the approximate problem in $\Gw^\ge=\Gw\setminus\overline B_\ge$,
\begin{equation}\label{glu1}
\left\{\BA {lll}
\CL_\gm u+g(u)=\gn_\gs&\quad\ {\rm in }\ \,   \Gw^\ge,\\[1mm]
\phantom{\CL_\gm +g(u)}
u=0&\quad\ {\rm on }\ \prt\Gw,\\[1mm]
\phantom{\CL_\gm +g(u)}
u=u_{k\gd_0}&\quad\ {\rm on }\ \prt B_\ge,
\EA\right.\end{equation}
where $u_{k\gd_0}$ is the solution of problem (\ref{sc1}) obtained in \rlemma{dirac}. It follows from \cite[Theorem 3.7]{V} that problem (\ref{glu1}) admits  a unique weak solution denoted by $U_{\gn_\gs,\ge}$, thanks to
 the fact that the operator is not singular in $\Gw^\ge$. We recall that $u_{\gn_\gs,\ge}$ is the solution of (\ref{sc10})
 and  $\BBG_{\gm}[\gd_0]$ the fundamental solution in $\Gw$. Then
\begin{equation}\label{glu2}\BA {lll}
\max\{u_{k\gd_0},u_{\gn_\gs,\ge}\}\leq U_{\gn_\gs,\ge}\leq u_{\gn_\gs}+k\BBG_{\gm}[\gd_0]\quad{\rm in }\,\;\Gw^\ge.
\EA\end{equation}
Furthermore one has $U_{\gn_\gs,\ge}\leq U_{\gn_\gs,\ge'}$ in $\Gw^\ge$, for $0<\ge'<\ge$. Since $u_{\gn_\gs}\leq u_{\gn}$ and
both $k\BBG_{\gm}[\gd_0]$ and $u_{\gn}$ belong to $L^1(\Gw,|x|^{-1}d\gg_\gm)$, then it follows by the monotone convergence theorem that $U_{\gn_\gs,\ge}$ converges in $L^1(\Gw,|x|^{-1}d\gg_\gm)$ and almost everywhere to some function $U_{\gn_\gs}\in L^1(\Gw,|x|^{-1}d\gg_\gm)$. Since $\Gg_\gm$ is a supersolution for equation
$\CL_\gm u+g(u)=0$ in $B_\gs$, for  $0<\ge_0<\gs$ there exists $c_8:=c_8(\ge_0,\gs)>0$ such that
$$u_{\gn_\gs}(x)\leq c_8|x|^{\gt_+}\quad{\rm for\, all }\,\;x\in B_{\ge_0}.$$
 For any $\gd>0$, there exists $\ge_0$ such that $u_{\gn_\gs}(x)\leq\gd\BBG_{\gm}[\gd_0](x)$ in $B_{\ge_0}$. Hence
 $u_{\gn_\gs}+k\BBG_{\gm}[\gd_0]\leq (k+\gd)\BBG_{\gm}[\gd_0]$ in $B_{\ge_0}$, which implies
\begin{equation}\label{glu3}
g(U_{\gn_\gs,\ge})\leq g((k+\gd)\BBG_{\gm}[\gd_0])
\quad{\rm in }\,\;B_{\ge_0}\setminus\overline B_{\ge}
\end{equation}
and
$$\BA {lll}\myint{\Gw}{}g((k+\gd)\BBG_{\gm}[\gd_0])d\gg_\gm(x)\leq \myint{B_1}{}g(\tfrac{k+\gd}{c_\gm}|x|^{\gt_-})|x|^{\gt_+}dx=|S^{N-1}|\myint{0}{1}g(\tfrac{k+\gd}{c_\gm}r^{\gt_-})r^{\gt_++N-1}dr\\[2mm]
\phantom{\myint{\Gw}{}g((k+\gd)\BBG_{\gm}[\gd_0])d\gg_\gm(x)}
=c_9\myint{\frac{k+\gd}{c_\gm}}{\infty}g(t)t^{-2+\frac{2}{\gt_-}}
=c_9\myint{\frac{k+\gd}{c_\gm}}{\infty}g(t)t^{-1-p^*_\gm}dt
\\[2mm]
\phantom{-----------\, }<\infty.
\EA$$
Now, using the local $\Gd_2$-condition,with $a'=\frac{k}{c_\gm}\ge_0^{\gt_-}$, we see that
\begin{equation}\label{glu4}g(U_{\gn_\gs,\ge})\leq g(u_{\gn_\gs}+\tfrac{k}{c_\gm}\ge_0^{\gt_-})\leq K(a')\left(g(u_{\gn_\gs})+g(a')\right)
\quad{\rm in }\,\;\Gw^{\ge_0}.
\end{equation}
From (\ref{glu3}) and (\ref{glu4}) we infer that $g(U_{\gn_\gs,\ge})$ is bounded in $L^1(\Gw^\ge,d\gg_\gm)$ independently of $\ge$. If
$\xi\in C_c^{1,1}(\Gw^*)$, we have for $\ge>0$ small enough so that supp$\,(\xi)\subset \Gw^\ge$
$$\myint{\Gw}{}\left(U_{\gn_\gs,\ge}\CL^*_\gm\xi+g(U_{\gn_\gs,\ge})\xi\right)d\gg_\gm=\myint{\Gw}{}\xi\Gg_\gm d\gn_\gs
$$
and letting $\ge\to 0$ we obtain that
\begin{equation}\label{glu5}
\myint{\Gw}{}\left(U_{\gn_\gs}\CL^*_\gm\xi+g(U_{\gn_\gs})\xi\right)d\gg_\gm=\myint{\Gw}{}\xi\Gg_\gm d\gn_\gs.
\end{equation}
Let $\xi\in C_c^{1,1}(\overline\Gw)$ and $\eta_n\in C^{1,1}(\BBR^N)$ a nonnegative cut-off function such that
$0\leq\eta_n\leq 1$, $\eta_n\equiv 1$ in $B^c_{\frac 2n}$, $\eta_n\equiv 0$ in $B_{\frac 1n}$, and choose $\xi\eta_n$ for test function. Then
\begin{equation}\label{glu6}
\myint{\Gw}{}\left(\eta_nU_{\gn_\gs}\CL^*_\gm\xi+g(U_{\gn_\gs})\eta_n\xi\right)d\gg_\gm -\myint{\Gw}{}U_{\gn_\gs}A_nd\gg_\gm=\myint{\Gw}{}\xi\eta_n\Gg_\gm d\gn_\gs
\end{equation}
with
\begin{equation}\label{glu7}
A_n=\xi\Gd\eta_n+2\langle\nabla\eta_n,\nabla\xi\rangle+2\gt_+\xi\langle\nabla\eta_n,\tfrac{x}{|x|^2}\rangle.
\end{equation}
Clearly
$$\lim_{n\to\infty}\myint{\Gw}{}\left(\eta_nU_{\gn_\gs}\CL^*_\gm\xi+g(U_{\gn_\gs})\eta_n\xi\right)d\gg_\gm=
\myint{\Gw}{}\left(U_{\gn_\gs}\CL^*_\gm\xi+g(U_{\gn_\gs})\xi\right)d\gg_\gm
$$
and
$$\lim_{n\to\infty}\myint{\Gw}{}\xi\eta_n\Gg_\gm d\gn_\gs=\myint{\Gw}{}\xi\Gg_\gm d\gn_\gs.
$$
We take
$$\eta_n(r)=\left\{\BA {lll}
\frac12-\frac12\cos\left(n\gp\left(r-\frac1n\right)\right)&\ \quad{\rm if}\,\;\frac1n\leq r\leq\frac2n,\\[1mm]
0&\ \quad{\rm if}\,\;r<\frac1n,\\[1mm]
1&\ \quad{\rm if}\,\;r>\frac2n.
\EA\right.
$$
Then
$$A_n=\frac{n^2\gp^2}{2}\cos\left(n\gp\left(r-\frac1n\right)\right)+\frac{n\gp}{2}\frac{N-1+2\gt_+}{r}\sin\left(n\gp\left(r-\frac1n\right)\right).
$$
Letting $\ge\to 0$ in (\ref{glu2}), we have
$$U_{\gn_\gs}(x)=k\BBG_{\gm}[\gd_0](x)(1+o(1))=\frac{k}{c_\gm}|x|^{\gt_-}(1+o(1))\qquad{\rm as}\,\; x\to 0.$$
Hence
\begin{equation}\label{glu8}\lim_{n\to\infty}\myint{\Gw}{}U_{\gn_\gs}A_nd\gg_\gm=\frac{2k|S^{N-1}|\sqrt{\gm-\gm_0}}{c_\gm}=k.
\end{equation}
This implies that $U_{\gn_\gs}$ is the solution of (\ref{eq1.1}) with $\gn$ replaced by $\gn_\gs+k\gd_0$.\qeda

\blemma{sigma0} Let $\gn=\gn\lfloor_{\Gw^*}+k\gd_0\in\overline{\frak M}_+(\Gw;\Gg_\gm)$. Then there exists a unique weak solution to (\ref{eq1.1}).
\es
\Proof Following the notations of \rlemma{sigma}, we set $\gn_\gs=\chi_{_{B_\gs}}\gn\lfloor_{\Gw^*}$ and denote by $U_{\gn_\gs}$ the solution of
\begin{equation}\label{glu9}
\left\{\BA {lll}
\CL_\gm u+g(u)=\gn_\gs+k\gd_0&\quad{\rm in }\,\ \Gw,\\[1.5mm]
\phantom{\CL_\gm +g(u)}
u=0&\quad{\rm on }\ \prt\Gw.
\EA\right.
\end{equation}
It is a positive function and there holds
\begin{equation}\label{glu10}\max\{u_{k\gd_0},u_{\gn_\gs}\}\leq U_{\gn_\gs}\leq u_{\gn_\gs}+k\BBG_{\gm}[\gd_0]\qquad{\rm in }\,\;\Gw.
\end{equation}
Since the mapping $\gs\mapsto U_{\gn_\gs}$ is decreasing, then there exists $U=\displaystyle \lim_{\gs\to 0}U_{\gn_\gs}$ and $U$ satisfies
(\ref{glu10}). As a consequence $U_{\gn_\gs}\to U$ in $L^1(\Gw,|x|^{-1}d\gg_\gm)$ as $\gs\to 0$. We take $\eta_1$ for test function in the weak formulation of (\ref{glu10}), then
$$\myint{\Gw}{}\left(|x|^{-1}U_{\gn_\gs}+\eta_1g(U_{\gn_\gs})\right) d\gg_\gm=\myint{\Gw}{}\eta_1\Gg_\gm d\gn_\gs+k\eta_1(0).
$$
By the monotone convergence theorem we obtain the identity
$$\myint{\Gw}{}\left(|x|^{-1}U+\eta_1g(U)\right) d\gg_\gm=\myint{\Gw}{}\eta_1d(\gg_\gm\gn\lfloor_{\Gw^*})+k\eta_1(0)=\myint{\Gw}{}\eta_1d(\gg_\gm\gn),
$$
and the fact that $g(U_{\gn_\gs})\to g(U)$ in $L^1(\Gw,\gr d\gg_\gm)$. Going to the limit as $\gs\to 0$ in the weak formulation of
(\ref{glu9}), we infer that $U=u_{\gn}$ is the solution of (\ref{eq1.1}).\qeda\medskip

\nind{\it Proof of Theorem B}. Assume $\gn=\gn\lfloor_{\Gw^*}+k\gd_0\in\overline{\frak M}(\Gw;\Gg_\gm)$ satisfies $k>0$ and let $\gn_+=\gn_+\lfloor_{\Gw^*}+k\gd_0$
and $\gn_-=\gn_-\lfloor_{\Gw^*}$ the positive and the negative part of $\gn$. We denote by $u_{\gn_+}$ and $u_{-\gn_-}$ the weak solutions of (\ref{eq1.1}) with respective data $\gn_+$ and $-\gn_-$. For $0<\ge<\gs$ such that $\overline B_\gs\subset\Gw$, we set
$\gn_\gs=\chi_{_{B_\gs}}\gn\lfloor_{\Gw^*}$, with positive and negative part $\gn_{\gs+}$ and $\gn_{\gs-}$ and denote by $U_{\gn_{\gs+},\ge}$, $U_{-\gn_{\gs-},\ge}$ and $U_{\gn_{\gs},\ge}$ the respective solutions of
\begin{equation}\label{B1}
\left\{\BA {lll}
\CL_\gm u+g(u)=\gn_{\gs+}&\quad{\rm in }\ \Gw^\ge,\\[1mm]
\phantom{\CL_\gm +g(u)}
u=0&\quad{\rm on }\ \prt\Gw,\\[1mm]
\phantom{\CL_\gm +g(u)}
u=u_{k\gd_0}&\quad{\rm on }\ \prt B_\ge,
\EA\right.\end{equation}

\begin{equation}\label{B2}\left\{\BA {lll}
\CL_\gm u+g(u)=-\gn_{\gs-}&\ \quad{\rm in }\,\ \Gw^\ge,\\[1mm]
\phantom{\CL_\gm +g(u)}
u=0&\ \quad{\rm on }\ \prt\Gw\cup \prt B_\ge
\EA\right.
\end{equation}
and
\begin{equation}\label{B3}\left\{\BA {lll}
\CL_\gm u+g(u)=\gn_{\gs}&\quad{\rm in }\,\ \Gw^\ge\\
\phantom{\CL_\gm +g(u)}
u=0&\quad{\rm on }\ \prt\Gw\\
\phantom{\CL_\gm +g(u)}
u=u_{k\gd_0}&\quad{\rm on }\ \prt B_\ge,
\EA\right.
\end{equation}
Then
\begin{equation}\label{B4}U_{-\gn_{\gs-},\ge}\leq U_{\gn_{\gs},\ge}\leq U_{\gn_{\gs+},\ge}.
\end{equation}
Furthermore $U_{\gn_{\gs+},\ge}$ satisfies (\ref{glu2}) and, in coherence with the notations of Lemma \ref{noatom} with $\gn_\gs$ replaced by $-\gn_{\gs-}$,
\begin{equation}\label{B5}
u_{-\gn_{\gs-}}\leq U_{-\gn_{\gs-},\ge}= u_{-\gn_{\gs-},\ge}.
\end{equation}
By compactness, $\{U_{\gn_{\gs},\ge_j}\}_{\ge_j}$ converges almost everywhere in $\Gw$ to some function $U$ for some sequence $\{\ge_j\}$ converging to $0$. Moreover $U_{\gn_{\gs},\ge_j}$ converges to $U_{\gn_\gs}$ in $L^1(\Gw,|x|^{-1}d\gg_\gm)$ because $U_{\gn_{\gs+},\ge}\to u_{\gn_{\gs+}+k\gd_0}$ and $u_{-\gn_{\gs-},\ge}\to u_{-\gn_{\gs-}}$ in $L^1(\Gw,|x|^{-1}d\gg_\gm)$ by  Lemma \ref{noatom} and (\ref{B4}) holds. Similarly $g(U_{\gn_{\gs},\ge_j})$ converges to $g(U)$ in $L^1(\Gw,\gr d\gg_\gm)$. This implies that $U$ satisfies
 $$\myint{\Gw}{}\left(U\CL^*_\gm\xi+g(U)\xi\right)d\gg_\gm=\myint{\Gw}{}\xi\Gg_\gm d\gn_\gs\quad{\rm for\ all}\ \,\xi\in C^{1,1}_c(\Gw^*).
 $$
 In order to use test functions in $C^{1,1}_c(\overline\Gw)$, we proceed as in the proof of \rlemma{sigma}, using the inequality (derived from (\ref{B4})) and the
 \begin{equation}\label{B6}
u_{-\gn_{\gs-}}\leq U_{\gn_\gs}\leq u_{\gn_{\gs+}+k\gd_0}.
\end{equation}
 By (\ref{glu4}), $u_{\gn_{\gs+}+k\gd_0}(x)=k\BBG_{\gm}[\gd_0](x)(1+o(1))$ when $x\to 0$ and $u_{-\gn_{\gs-}}=o(\BBG_{\gm}[\gd_0])$ near $0$. This implies
 $U_{\gn_\gs}(x)=k\BBG_{\gm}[\gd_0](x)(1+o(1))$ as $x\to 0$ and we conclude as in the proof of \rlemma{sigma} that $u=u_{\gn_\gs+k\gd_0}$.\\
At end we let $\gs\to 0$. Up to a sequence $\{\gs_j\}$ converging to $0$ such that $u_{\gn_{\gs_j}+k\gd_0}\to U$ almost everywhere and
 \begin{equation}\label{B7}
u_{-\gn_{-}}\leq U\leq u_{\gn_{+}+k\gd_0}.
\end{equation}
Since by \rlemma{sigma0}, $u_{\gn_{\gs+}+k\gd_0}\to u_{\gn_{+}+k\gd_0}$ in $L^1(\Gw,|x|^{-1}d\gg_\gm)$ and $g(u_{\gn_{\gs+}+k\gd_0})\to g(u_{\gn_{+}+k\gd_0})$ in $L^1(\Gw,\gr d\gg_\gm)$, we infer that the convergences of $u_{\gn_{\gs_j}+k\gd_0}\to U$ and
$g(u_{\gn_{\gs_j}+k\gd_0})\to g(U)$ occur respectively in the same space, therefore $U=u_{\gn+k\gd_0}$, it is the weak solution of
(\ref{eq1.1}).\qeda\\

\nind\Remark In the course of the proof we have used the following result which is independent of any assumption on $g$ but for the monotonicity: If $\{\gn_n\}\subset \overline{\mathfrak M}_+(\Gw;\Gg_\gm)$ is an increasing sequence of $g$-good measures converging to a measure $\gn\in \overline{\mathfrak M}_+(\Gw;\Gg_\gm)$, then $\gn$ is a  $g$-good measure, $\{u_{\gn_n}\}$ converges to
$u_{\gn}$ in $L^1(\Gw,|x|^{-1}d\gg_\gm)$ and $\{g(u_{\gn_n})\}$ converges to
$g(u_{\gn})$ in $L^1(\Gw,\gr d\gg_\gm)$.


\subsection{Proof of Theorem C}
The construction of a solution is essentially similar to the one of Theorem B, the only modifications lies in \rlemma{sigma}. Estimate (\ref{glu2}) remains valid with
 \begin{equation}\label{B8}
u_{k\gd_0}(x)=\tfrac{k}{|S^{N-1}|}|x|^{\frac{2-N}{2}}\ln |x|^{-1}(1+o(1))=k\BBG_{\gm}[\gd_0](x)(1+o(1))\quad{\rm as}\,\;x\to 0.
\end{equation}
Since $u_{\gn_\gs}(x)\leq c|x|^{\frac{2-N}{2}}$, (\ref{glu3}) holds with $\gd>0$ arbitrarily small. Next
$$\BA {lll}\myint{\Gw}{}g((k+\gd)\BBG_{\gm}[\gd_0])d\gg_\gm(x)\leq \myint{B_1}{}g\left(\tfrac{k+\gd}{|S^{N-1}|}|x|^{\frac{2-N}{2}}\ln |x|^{-1}\right)|x|^{\frac{2-N}{2}}dx\\[2mm]\phantom{\myint{\Gw}{}g((k+\gd)\BBG_{\gm}[\gd_0])d\gg_\gm(x)}
=|S^{N-1}|\myint{0}{1}g\left(\tfrac{k+\gd}{|S^{N-1}|}r^{\frac{2-N}{2}}\ln r^{-1}\right)r^{\frac N2}dr
\\[2mm]\phantom{\myint{\Gw}{}g((k+\gd)\BBG_{\gm}[\gd_0])d\gg_\gm(x)}
=c_{10}\myint{c'}{\infty}g(t\ln t)t^{-\frac{2N}{N-2}}<\infty,
\EA$$
by (\ref{sc6-2}) and (\ref{1.4-3}).  The end of the proof for Theorem C is similar to the one of Theorem B.\qeda\medskip

\nind{\it Proof of Corollary D}. If $g(r)=g_p(r)=|r|^{p-1}r$, $p>1$, the existence of a solution with $\gn=k\gd_0$ is obtained if
$$\myint{1}{\infty}(t^p-|t|^p)t^{p_{\gm_*}}dt<\infty\;\,{\rm if }\;\gm>\gm_0\;\,{\rm and }\;\,
\myint{1}{\infty}(t^p-|t|^p)t^{-\frac{2N}{N-2}}(\ln t)^{\frac{N+2}{N-2}}dt<\infty\;\,{\rm if }\;\gm=\gm_0.
$$
In both case $p<p_{\gm_*}$. If $k=0$ and $\gn\lfloor_{\Gw^*}\neq 0$, the existence is ensured if (\ref{eq1.1-2}) holds, hence $p<\frac{N}{N-2}$. Assertion (iii) follows.\qeda
\setcounter{equation}{0}
\section{The supercritical case}
\subsection{Reduced measures}
The notion of reduced measures introduced by Brezis, Marcus and Ponce \cite{BMP} turned out to be a useful tool in the construction of solutions in a measure framework. We will develop only the aspect needed for the proof of theorem E. If $k\in\BBN$, we set
\begin{equation}\label{red1}
g_k(r)=\left\{\BA{lll}\min\{g(r),g(k)\}&\quad{\rm if}\,\; r\geq 0,\\[1.5mm]
\max\{g(r),g(-k)\}&\quad{\rm if}\,\; r>0.
\EA\right.
\end{equation}
Since $g_k$ satisfies (\ref{1.4-2}) and (\ref{1.4-3}), for any $\gn\in \overline{\mathfrak M}_+(\Gw;\Gg_\gm)$ there exists a unique weak solution $u=u_{\gn,k}$ of
\begin{equation}\label{red2}
\left\{\BA {lll}
\CL_\gm u+g_k(u)=\gn&\quad{\rm in}\,\ \Gw,\\[1.5mm]
\phantom{\CL_\gm+g_k(u)}
u=0&  \quad{\rm on} \ \prt\Gw.
\EA\right.\end{equation}
Furthermore,   from the proof of \rlemma{sigma0} and Kato's type estimates Proposition \ref{pr 2.1} we have that
\begin{equation}\label{red3}\BA {lll}
0\le u_{\gn_+,k'}\leq u_{\gn_+,k}&\ \quad{\rm for\, all}\,\ k'\geq k>0.
\EA\end{equation}


\bprop {redprop}Let $\gn\in \overline{\mathfrak M}_+(\Gw;\Gg_\gm)$. Then the sequence of weak solutions $\{u_{\gn,k}\}$ of
\begin{equation}\label{red4}
\left\{\BA {lll}
\CL_\gm u+g_k(u)=\gn&\quad{\rm in}\,\ \Gw,\\[1.5mm]
\phantom{\CL_\gm+g_k(u)}
u=0&\quad{\rm on}\ \prt\Gw
\EA\right.\end{equation}
decreases and converges, when $k\to\infty$, to some nonnegative function $u$ and there exists a measure $\gn^*\in \overline{\mathfrak M}_+(\Gw;\Gg_\gm)$ such that $0\leq\gn^*\leq\gn$ and $u=u_{\gn^*}$.
\es
\Proof The proof is similar to the one of \cite[Prop. 4.1]{BMP}.  Observe that
$u_{\gn,k}\downarrow u^*$ and the sequence $\{u_{\gn,k}\}$ is uniformly integrable in $L^1(\Gw,|x|^{-1} d\gg_\gm)$. By Fatou's lemma
$u$ satisfies
\begin{equation}\label{red5}
\myint{\Gw}{}\left(u^*\CL_{\gm}^*\xi+g(u^*)\xi\right) d\gg_\gm(x)\leq \myint{\Gw}{}\xi d(\Gg_\gm \gn)
\quad{\rm for\, all}\,\; \xi\in\BBX\gm(\Gw),\;\xi\geq 0.
\end{equation}
Hence $u^*$ is a subsolution of (\ref{eq1.1}) and by construction it is the largest of all nonnegative subsolutions. The mapping
$$\xi\mapsto \myint{\Gw}{}\left(u^*\CL_{\gm}^*\xi+g(u^*)\xi\right) d\gg_\gm(x)\qquad{\rm for\, all}\,\;\xi\in C^{\infty}_c(\Gw),
$$
is a positive distribution, hence a measure $\nu^*$, called {\it the reduced measure of $\gn$}.  It satisfies $0\leq\gn^*\leq\gn$ and
$u^*=u_{\nu^*}$. \qeda

\blemma{redlem2} Let $\gn,\gn'\in \overline{\mathfrak M}_+(\Gw;\Gg_\gm)$. If $\gn'\leq \gn$ and $\gn=\gn^*$, then $\gn'=\gn'^*$.
\es
\Proof Let $u_{\gn',k}$ be the weak solution of the truncated equation
 \begin{equation}\label{red6}\left\{\BA {lll}
\CL_\gm u+g_k(u)=\gn'&\quad{\rm in}\,\ \Gw,\\[1.5mm]
\phantom{\CL_\gm+g_k(u)}
u=0&\quad{\rm on}\ \prt\Gw.
\EA\right.
\end{equation}
 Then $0\leq u_{\gn',k}\leq u_{\gn,k}$. By \rprop{redprop}, we know that $u_{\gn,k}\downarrow u_{\gn^*}=u_\gn$ and $u_{\gn',k}\downarrow u'^*$ a.e. in
 $L^1(\Gw,|x|^{-1} d\gg_\gm)$ and then
 $$\CL_\gm(u_{\gn,k}-u_\gn)+g_k(u_{\gn,k})-g_k(u_\gn)=g(u_\gn)-g_k(u_\gn),
 $$
 hence, by Proposition \ref{pr 2.1},
 $$\myint{\Gw}{}(u_{\gn,k}-u_\gn))|x|^{-1}d\gg_\gm+\myint{\Gw}{}|g_k(u_{\gn,k})-g_k(u_\gn)|\eta_1d\gg_\gm\leq
 \myint{\Gw}{}|g(u_\gn)-g_k(u_\gn)|\eta_1d\gg_\gm.
 $$
 By the increasing monotonicity of mapping  $k\mapsto g_k(u_\gn)$,  we have $g_k(u_\gn)\to g(u_\gn)$ in $L^1(\Gw,\gr d\gg_\gm)$ as $k\to+\infty$, hence
 $$\myint{\Gw}{}|g_k(u_{\gn,k})-g(u_\gn)|\eta_1d\gg_\gm\leq 2\myint{\Gw}{}|g(u_\gn)-g_k(u_\gn)|\eta_1d\gg_\gm\to 0\;\;{\rm as}\;\,n\to\infty.
 $$
 Because $g_k(u'_{\gn,k})\leq g_k(u_{\gn,k})$ it follows by Vitali's convergence theorem that $g_k(u_{\gn',k})\to g(u'^*)$ in $L^1(\Gw,\gr d\gg_\gm)$. Using the weak formulation of  (\ref{red6}), we infer that $u'^*$ verifies
$$\myint{\Gw}{}\left(u'^*\CL_\gm^*\xi +g(u'^*)\xi\right)d\gg_\gm=\myint{\Gw}{}\xi d(\gg_\gm\gn')\quad\ {\rm for\ all\ }  \xi\in\BBX_\gm(\Gw).
$$
 This yields $u'^*=u_{\gn'}$. \qeda\medskip

 The next result follows from \rlemma{redlem2}.

\blemma{redlem3} Assume that $\gn=\gn\lfloor_{\Gw^*}+k\gd_0\in \overline{\mathfrak M}_+(\Gw;\Gg_\gm)$, then
$\gn^*=\gn^*\lfloor_{\Gw^*}+k^*\gd_0\in \overline{\mathfrak M}_+(\Gw;\Gg_\gm)$ with $\gn^*\lfloor_{\Gw^*}\leq \gn\lfloor_{\Gw^*}$
and $k^*\leq k$. More precisely, \smallskip\smallskip

\nind (i) If $\gm>\gm_0$ and $g$ satisfies (\ref{1.4-2}), then $k=k^*$.\\[1mm]
\nind (ii) If $\gm=\gm_0$ and $g$ satisfies (\ref{1.4-3}), then $k=k^*$.\\[1mm]
\nind (ii) If $\gm>\gm_0$ (resp. $\gm=\gm_0$)  and $g$ does not satisfy (\ref{I-2}) (resp. (\ref{1.4-3})), then $k^*=0$.
\es
\medskip

The next result is useful in applications.
\begin{corollary} \label{redlarge} If $\gn\in \overline{\mathfrak M}_+(\Gw;\Gg_\gm)$, then $\gn^*$ is the largest $g$-good measure  smaller or equal to $\gn$.
\end{corollary}
\Proof Let $\gl\in \overline{\mathfrak M}_+(\Gw;\Gg_\gm)$ be a $g$-good measure, $\gl\leq\gn$. Then $\gl^*=\gl\leq\gn^*$. Since $\gn^*$ is a $g$-good measure the result follows.\qeda

\medskip

\nind{\it Proof of Theorem E}. Assume that $\gn\geq 0$. By \rlemma {redlem2} and Remark at the end of Section 3.5 the following assertions are equivalent:\\
(i) $\gn$ is $g_p$-good.\\
(ii) For any $\gs>0$, $\gn_\gs=\chi_{_{B^c_\gs}}\gn$ is $g_p$-good. \smallskip

 If $\gn_\gs$ is good, then $u_{\gn_\gs}$ satisfies
\bel{E1}-\Gd u_{\gn_\gs}+u^p_{\gn_\gs}=\gn_\gs-\myfrac{\gm}{|x|^2}u_{\gn_\gs}\quad {\rm in}\ \CD'(\Gw^*)
\ee
  and since $u_{\gn_\gs}(x)\leq c|x|^{\gt_+}$ if $|x|\leq \frac{\gs}{2}$
(\ref{E1}) holds in $\CD'(\Gw)$.  This implies that $u\in L^p(\Gw)$  and $|x|^{-2}u_{\gn_\gs}\in L^\ga(B_{\frac{\gs}{2}})$ for any
$\ga<\frac{N}{(2-\gt_+)_+}$. Using \cite{BP2} the measure $\gn_\gs$ is absolutely continuous with respect to the $c_{2,p'}$-Bessel capacity. If $E\subset\Gw$ is a Borel set such that $c_{2,p'}(E)=0$, then $c_{2,p'}(E\cap B^c_\gs)=0$, hence $\gn (E\cap B^c_\gs)=\gn_\gs(E\cap B^c_\gs)=0$. By the monotone convergence theorem $\gn(E)=0$. \smallskip

 Conversely, if $\gn$ is nonnegative and absolutely continuous with respect to the $c_{2,p'}$-Bessel capacity, then so is $\gn_\gs =\chi_{_{B^c_\gs}}\gn$. For $0\leq\ge\leq\frac\gs 2$ we consider the problem
\begin{equation}\label{F1}
\left\{\BA {lll}
\,-\Gd u+\myfrac{\gm}{|x|^2}u+u^p=\gn_\gs&\ \ {\rm in\,\ }\Gw^\ge:=\Gw\setminus B_\ge,\\[1.5mm]
\phantom{--------}u=0&\ \ {\rm on\,\ }\prt  B_\ge,\\[1mm]
\phantom{--------}u=0&\ \ {\rm on\,\ }\prt  \Gw.
\EA\right.\end{equation}
Since $\tfrac{\gm}{|x|^2}$ is bounded in $\Gw^\ge$ and $\gn_\gs$  is absolutely continuous with respect to the $c_{2,p'}$ capacity there exists a solution $u_{\gn_\gs,\ge}$ thanks to \cite {BP2}, unique by monotonicity. Now the mapping $\ge\mapsto u_{\gn_\gs,\ge}$ is decreasing. We use the method developed in Lemma \ref{noatom}, when $\ge\to 0$, we know that $u_{\gn_\gs,\ge}$ increase to some $u_\gs$ which is dominated by $\BBG[\gn_\gs]$ and satisfies
\bel{F2}\left\{\BA {lll}
-\Gd u+\myfrac{\gm}{|x|^2}u+u^p=\gn_\gs&\ \ {\rm in\,\ }\Gw^*,\\[1.5mm]
\phantom{--------}u=0&\ \ {\rm on\ }\prt  \Gw.
\EA\right.\ee
Because $u_\gs\leq \BBG[\gn_\gs]$ and $\gn_\gs=0$ in $B_\gs$, there holds $u(x)\leq c'_{11}\Gg_\gm(x)$ in $B_{\frac\gs 2}$, and then
$u_\gs$ is a solution in $\Gw$ and $u=u_{\gn_\gs}$. Letting $\gs\to 0$, we conclude as in Lemma \ref{noatom}  that $u_{\gn_\gs}$ converges to $u_\gn$ which is the weak solution of
\bel{F3}\left\{\BA {lll}
-\Gd u+\myfrac{\gm}{|x|^2}u+u^p=\gn&\ \ {\rm in\,\ }\Gw,\\[1.5mm]
\phantom{--------}u=0&\ \ {\rm on\ }\prt  \Gw.
\EA\right.\ee
If $\gn$ is a signed measure absolutely continuous with respect to the $c_{2,p'}$-capacity, so are $\gn_+$ and $\gn_-$. Hence there exists solutions $u_{\gn_+}$ and $u_{\gn_-}$. For $0<\ge<\frac\gs 2$ we construct $u_{\gn_\gs,\ge}$ with the property that
$-u_{-\gn_{-\,\gs},\ge}\leq u_{\gn_\gs,\ge}\leq u_{\gn_{+\,\gs},\ge}$, we let $\ge\to 0$ and derive the existence of $u_{\gn_\gs}$ which is eventually the weak solution of
\bel{F4}\left\{\BA {lll}
-\Gd u+\myfrac{\gm}{|x|^2}u+|u|^{p-1}u=\gn_\gs&\ \ {\rm in\,\ }\Gw^*,\\[1.5mm]
\phantom{----------}u=0&\ \ {\rm on\ }\prt  \Gw,
\EA\right.\ee
and satisfies $-u_{-\gn_{-\,\gs}}\leq u_{\gn_\gs}\leq u_{\gn_{+\,\gs}}$. Letting $\gs\to 0$ we then derive that $\displaystyle u=\lim_{\gs\to 0}u_{\gn_\gs}$ satisfies
\bel{F5}\left\{\BA {lll}
-\Gd u+\myfrac{\gm}{|x|^2}u+|u|^{p-1}u=\gn&\ \ {\rm in\,\ }\Gw^*,\\[1.5mm]
\phantom{----------}u=0&\ \ {\rm on\ }\prt  \Gw.
\EA\right.\ee
Hence $u=u_\gn$ and $\gn$ is a good solution.\qeda
\medskip

\nind{\it Proof of Theorem F}.  {\bf Part 1.}  Without loss of generality we can assume that $\Gw$ is a bounded smooth domain. Let $K\subset\Gw$ be compact. If $0\in K$ and $p< p^*_{\gm}$ there exists a solution $u_{k\gd_0}$, hence
$K$ is not removable. If $0\notin K$ and $c_{2,p'}(K)>0$, there exists a capacitary measure $\gn_K\in W^{-2,p}(\Gw)\cap \mathfrak M_+(\Gw)$ with support in $K$. This measure is $g_p$-good by Theorem E, hence $K$ is not removable. \smallskip

 {\bf Part 2.}   Conversely we first assume that $0\notin K$. Then there exists a subdomain $D\subset \Gw$ such that $0\notin \bar D$ and $K\subset D$. Hence a solution $u$ of (\ref{1p.2}) is also a solution of
$$-\Gd u+\frac{\gm}{|x|^2}u+|u|^{p-1}u=0\ \quad{\rm in\,\ }D\setminus K
$$
and the coefficient $\frac{\gm}{|x|^2}$ is uniformly bounded in $\bar D$. By \cite[Theorem 3.1]{BP2} it can be extended as a $C^2$ solution of the same equation in $\Gw'$. Hence $K$ is removable if $c_{2,p'}(K)=0$.

 If $0\in K$ we have to assume at least $p\geq p^*_{\gm}$ in order  that $0$ is removable and $p\geq p_0$ in order there exists non-empty set with zero $c_{2,p'}$-capacity. Let $\gz\in C_c^{1,1}(\Gw)$ with $0\leq\gz\leq 1$, vanishing in a compact neighborhood $D$ of $K$. Since $0\notin \Gw\setminus D$, we first consider the case where $u$ is nonnegative and satisfies in the usual sense
$$-\Gd u+\frac{\gm}{|x|^2}u+ u^p=0\quad{\rm in\,\;}\Gw\setminus D.
$$
Taking $\gz^{2p'}$ for test function, we get
$$-2p'\myint{\Gw}{}u\gz^{2p'-1}\Gd\gz dx-2p'(2p'-1)\myint{\Gw}{}u\gz^{2p'-2}|\nabla\gz|^2dx+\gm\myint{\Gw}{}
\frac{u\gz^{2p'}}{|x|^2}dx+\myint{\Gw}{}\gz^{2p'}u^p dx=0.
$$
There holds
$$\left|\myint{\Gw}{}u\gz^{2p'-1}\Gd\gz dx\right|\leq \left(\myint{\Gw}{}\gz^{2p'}u^p dx\right)^{\frac 1p}
\left(\myint{\Gw}{}|\Gd\gz|^{p'}\gz^{p'}dx\right)^{\frac 1{p'}},
$$
$$0\leq \myint{\Gw}{}u\gz^{2p'-2}|\nabla\gz|^2dx\leq \left(\myint{\Gw}{}\gz^{2p'}u^p dx\right)^{\frac 1p}
\left(\myint{\Gw}{}|\nabla\gz|^{2p'}dx\right)^{\frac 1{p'}},
$$
and
$$0\leq \myint{\Gw}{}
\frac{u\gz^{2p'}}{|x|^2}dx\leq \left(\myint{\Gw}{}\gz^{2p'}u^p dx\right)^{\frac 1p}
\left(\myint{\Gw}{}\myfrac{\gz^{2p'}}{|x|^{2p'}}dx\right)^{\frac 1{p'}}.
$$
By standard estimates and Gagliardo-Nirenberg inequality (and since $0\leq\gz\leq 1$),
$$\left(\myint{\Gw}{}|\Gd\gz|^{p'}\gz^{p'}\right)^{\frac 1{p'}}\leq c_{11} \norm\gz_{W^{2,p'}}
$$
and
$$\left(\myint{\Gw}{}|\nabla\gz|^{2p'}dx\right)^{\frac 1{p'}}\leq c_{12} \norm\gz_{W^{2,p'}}.
$$
Finally, if $p>p_0:=\frac{N}{N-2}$, then $2p'<N$ which implies that there exists $c_{13}$ independent of $\gz$ (with value in $[0,1]$) such that
$$\left(\myint{\Gw}{}\myfrac{\gz^{2p'}}{|x|^{2p'}}dx\right)^{\frac 1{p'}}\leq \left(\myint{B_1}{}\myfrac{dx}{|x|^{2p'}}\right)^{\frac 1{p'}}:=c_{13}.
$$
Next we set
$$X=\left(\myint{\Gw}{}\gz^{2p'}u^p dx\right)^{\frac 1p}
$$
and we obtain if $\gm\geq 0$, if $p\geq p_0$
\bel{E2}
X^p-\left(2p'(2p'-1)c_{12}-p'c_{12}\right)\norm\gz_{W^{2,p'}}X\leq 0
\ee
and if $\gm<0$ if $p> p_0$
\bel{E3}
X^p-\left(\left(2p'(2p'-1)c_{12}-p'c_{12}\right)\norm\gz_{W^{2,p'}}-c_{13}\gm\right)X\leq 0.
\ee
However, the condition $p>p_0$ is ensured when $\gm<0$ since $p\geq p_\gm^*>p_0$.
We consider a sequence $\{\eta_n\}\subset \CS(\BBR^N)$ such that $0\leq\eta_n\leq 1$, $\eta_n=0$ on a neighborhood of $K$ and such that $\norm{\eta_n}_{W^{2,p'}}\to 0$ when $n\to\infty$. Such a sequence exists by \cite{Me} result since $c_{2,p'}(K)=0$. Let $\xi\in C^{\infty}_0(\Gw)$ such that $0\leq \xi\leq 1$ and with value $1$ in a neighborhood of $K$. We take $\gz:=\gz_n=(1-\eta_n)\xi$ in the above estimates. Letting $n\to\infty$, then $\gz_n\to \xi$ in $W^{2,p'}$ and finally
\bel{E4}X^{p-1}=\left(\myint{\Gw}{}\xi^{2p'}u^p dx\right)^{\frac {p-1}p}\leq \left(2p'(2p'-1)c_{12}-p'c_{12}\right)\norm\xi_{W^{2,p'}}+c_{13}\gm_-;
\ee
under the condition that $p>p_0$ if $\gm<0$, in which case there also holds
\bel{E5}
\myint{\Gw}{}
\frac{u\gz^{2p'}}{|x|^2}dx\leq c_{13}X.
\ee
However the condition $p>p_0$ is not necessary in order the left-hand side of (\ref{E5}) be bounded, since we have
\bel{E6}
\gm\myint{\Gw}{}
\frac{u\gz^{2p'}}{|x|^2}dx+X^p\leq \left(2p'(2p'-1)c_{12}-p'c_{12}\right)\norm\gz_{W^{2,p'}}X,
\ee
and $X$ is bounded.

Next we take $\gz:=\gz_n=(1-\eta_n)\xi$ for test function in (\ref{1p.2}) and get
$$
-\myint{\Gw}{}\left( (1-\eta_n)\Gd\xi-\xi\Gd\eta_n-2\langle\nabla\eta_n,\nabla\xi\rangle \right) udx+
\gm\myint{\Gw}{}
\frac{u\gz_n}{|x|^2}dx+\myint{\Gw}{}\gz_nu^p dx=0.
$$
Since
$$\myint{\Gw}{}u\xi\Gd\eta_n dx\leq \left(\myint{\Gw}{}u^p\xi dx\right)^{\frac1p}\norm{\eta_n}_{W^{2,p'}}\to 0\quad{\rm as\;\,}n\to\infty
$$
and
$$\left|\myint{\Gw}{}u\langle\nabla\eta_n,\nabla\xi\rangle dx\right|\leq \left(\myint{\Gw}{}u^p|\nabla\xi| dx\right)^{\frac1p}
\norm{\nabla\xi}_{L^\infty}\norm{\eta_n}_{W^{1,p'}}\quad{\rm as\;\,}n\to\infty,
$$
then we conclude that $u$ satisfies
\bel{E7}
-\myint{\Gw}{}u\Gd\xi dx+\gm\myint{\Gw}{}\frac{u\xi}{|x|^2}dx+\myint{\Gw}{}\xi u^p dx=0,
\ee
which proves that $u$ satisfies the equation in the sense of distributions. By standard regularity $u$ is $C^2$ in $\Gw^*$, and by the maximum principle $u(x)\leq c_{14}\Gg_\gm(x)$ in $B_{r_0}\subset\Gw$. Integrating by part as in the proof of \rlemma{dirac} we obtain that $u$ satisfies
\bel{E8}
\myint{\Gw}{}\left(u\CL^*_\gm\xi +\xi u^p \right)d\gg_\gm(x)=0\quad{\rm for\ every\ } \xi\in\BBX_\gm(\Gw).
\ee

Finally, if $u$ is a signed solution, then $|u|$ is a subsolution. For $\ge>0$ we set $K_\ge=\{x\in \BBR^N:\dist(x,K)\leq\ge\}$. If $\ge$ is small enough $K_\ge\subset\Gw$. Let $v:=v_\ge$ be the solution of
\bel{E9}\left\{\BA {lll}
-\Gd v+\myfrac{\gm}{|x|^2}v+v^p=0&\quad{\rm in\ }\Gw\setminus K_\ge,\\[1.5mm]
\phantom{-------\ \,}v=|u|\lfloor_{\prt K_\ge}&\quad{\rm on\ }\prt  K_\ge,\\[1.5mm]
\phantom{-------\ \,}v=|u|\lfloor_{\prt \Gw}&\quad{\rm on\ }\prt  \Gw.
\EA\right.\ee
Then $|u|\leq v_\ge$. Furthermore, by Keller-Osserman estimate as in \cite{GV}, there holds
\bel{E10}\BA {lll}
v_\ge(x)\leq c_{15}\dist (x,K_\ge)^{-\frac{2}{p-1}}\quad{\rm for\,all\, \,}x\in\Gw\setminus K_\ge,
\EA\ee
where $c_{14}>0$ depends on $N$, $p$ and $\gm$. Using local regularity theory and the Arzela-Ascoli theorem, there exists a sequence
$\{\ge_n\}$ converging to $0$ an a function $v\in C^2(\Gw\setminus K)\cap C(\bar\Gw\setminus K)$ such that $\{v_{\ge_n}\}$ converges to $v$ locally uniformly in $\bar\Gw\setminus K$ and in the $C^2_{loc}(\Gw\setminus K)$-topology. This implies that
$v$ is a positive solution of (\ref{1p.2}) in $\Gw\setminus K$. Hence it is a solution in $\Gw$. This implies that $u\in L^p(\Gw)$ and
$ |u(x)|\leq v(x)\leq c_{14}\Gg_\gm(x)$ in $\Gw^*$. We conclude as in the nonnegative case that $u$ is a weak solution in $\Gw$.\qeda
\bigskip

  \noindent{\bf Acknowledgements:}  H. Chen is supported by
NSF of China, No: 11726614, 11661045, by the Jiangxi Provincial Natural Science Foundation, No: 20161ACB20007, 
and by the Alexander von Humboldt Foundation.

\end{document}